\documentclass[final,onefignum,onetabnum,dvipsnames]{siamart171218}
\usepackage[utf8]{inputenc}
\usepackage{amsmath}
\usepackage{amsfonts}
\usepackage{amssymb}
\usepackage{amsopn}
\usepackage{graphicx}
\usepackage{bm}
\usepackage{enumitem}
\usepackage[font=footnotesize]{caption}
\usepackage{graphicx}
\usepackage{xcolor}
\usepackage{placeins}
\usepackage{ifthen}
\usepackage{xifthen}
\usepackage{subfig}

\usepackage{epstopdf}
\ifpdf
  \DeclareGraphicsExtensions{.eps,.pdf,.png,.jpg}
\else
  \DeclareGraphicsExtensions{.eps}
\fi

\newsiamremark{remark}{Remark}
\newsiamremark{hypothesis}{Hypothesis}
\crefname{hypothesis}{Hypothesis}{Hypotheses}
\newsiamthm{claim}{Claim}
\newsiamremark{assumption}{Assumption} 
\crefname{assumption}{Assumption}{Assumption} 

\headers{Path-following methods for calculating dispersion relations}{Peter Maxwell and Simen \AA.\ Ellingsen}

\title{Path-following methods for calculating linear surface wave dispersion relations on vertical shear flows.\thanks{arXiv PREPRINT VERSION \funding{This work was funded by the Norwegian Research Council (FRINATEK) {\#}249740.}}}

\author{Peter Maxwell\thanks{Norwegian University of Science and Technology, Department of Energy and Process Engineering (\email{peter.maxwell@ntnu.no}, \email{simen.a.ellingsen@ntnu.no}; \url{https://www.ntnu.edu/employees/peter.maxwell}, \url{https://www.ntnu.edu/employees/simen.a.ellingsen}).}
\and
Simen {\AA}.\ Ellingsen\footnote[2].}

\graphicspath{  {./figures/} }

\RequirePackage{waves_path_follow}

\newlist{pfgsteplist}{enumerate}{1}
\setlist[pfgsteplist]{label={\roman*.}}
\crefname{pfgsteplisti}{step}{steps}

\begin{document}

\maketitle

\begin{abstract}
The path-following scheme in Loisel and Maxwell [\textit{SIAM J. Matrix Anal. Appl.} 39-4 (2018), pp.\ 1726-1749] is adapted to efficiently calculate the dispersion relation curve for linear surface waves on an arbitrary vertical shear current. This is equivalent to solving the Rayleigh instability equation with linearised free-surface boundary condition for each sought point on the curve. Taking advantage of the analyticity of the dispersion relation, a path-following or continuation approach is adopted. The problem is discretized using a collocation scheme, parametrised along either a radial or angular path in the wave vector plane, and differentiated to yield a system of ODEs. After an initial eigenproblem solve using QZ decomposition, numerical integration proceeds along the curve using linear solves as the Runge--Kutta $F(\cdot)$ function; thus, many QZ decompositions on a $2N$ companion matrix are exchanged for one QZ decomposition and a small number of linear solves on a size $N$ matrix. A piecewise interpolant provides dense output. The integration represents a nominal setup cost afterwhich very many points can be computed at negligible cost whilst preserving high accuracy. Furthermore, a 2-dimensional interpolant suitable for scattered data query points in the wave vector plane is described. Finally, a comparison is made with existing numerical methods for this problem, revealing that the path-following scheme is asymptotically two orders of magnitude faster in number of query points.
\end{abstract}



\section{Introduction}
\label{sec_intro}
Engineering and the natural sciences are replete with eigenproblems for ordinary differential operators which depend on a finite set of parameters. We are interested in problems which are parametrised by a single real variable.

The canonical solution approach involves conversion to an algebraic problem via spatial discretization, which often leads to polynomial or even nonlinear eigenproblems of potentially large dimension. These can be solved using classical techniques for each sought parameter value. This strategy may become prohibitively expensive when the computation must be repeated many times. It can also be difficult to take advantage of the nearness of solutions for small parameter variations, forcing full calculations for each point in the parameter space. An alternative approach is to solve the discretized eigenproblem once for a fixed parameter value then use the local piecewise analyticity of the eigenvalue and eigenvector \cite{Kato1995}, \cite{AndrewChuLancaster1993}, \cite{QianChuTan2015} to calculate along the solution curve using a path-following or numerical continuation algorithm.

In a more general setting, this comprises a numerical continuation method whereby the parameter-dependent solution is calculated as an implicitly defined curve \cite{AllgowerGeorg2003}. Homotopy methods have a similar philosphy but introduce an artificial parameter to parametrise a convex homotopy to map from the solution of an `easy' problem to the solution of the actual problem \cite{Liao2003}. These methods tend to use predictor-corrector schemes such as pseudo-arclength continuation or similar approaches. We make reference to the homotopy method in \cite{LuiKellerKwok1997} and the invariant subspace methods in \cite{BeynThummler2010}, \cite{BeynEffenbergerKressner2011} as relevant examples. For a recent approach that shares a strong philosophical similarity with the material herein for working with time-varying matrix eigenproblems, albeit using different techniques (look-ahead finite difference formulas), see \cite{UhligYunong2019} and \cite{Uhlig2019}.

This paper is concerned with repurposing the path-following technique used in \cite{LoiselMaxwell2018} to solve a specific classical problem from wave-current interactions \cite{Rayleigh1979}, \cite{Rayleigh1887}, \cite{Rayleigh1895}, \cite{VelthuizenWijngaarden1969}, \cite{Yih1972}, \cite[s.\ IV]{Peregrine1976}, \cite{Shrira1993}, \cite{Zhang2005}, \cite{Miles2001a}, \cite{Miles2001b} : that of calculating the dispersion relation for perturbative linear order free-surface waves travelling atop a vertical shear flow. The problem is particularly suited as a motivating example of the technique: it is conceptually simple, it has an eigenvalue-dependent boundary condition, it is well-known from both the waves literature and hydrodynamic stability, and there is a practical requirement for efficient numerical solution.

We summarise our approach as follows. The original eigenproblem is spatially discretized using a collocation method, implicitly incorporating the boundary conditions, to obtain a parameter-varying system of equations that are then differentiated to yield an under-determined system of ODEs. An additional constraint is then included. After performing an initial eigenproblem solve, numerical integration can proceed along the solution curve using linear solves as the Runge--Kutta $F(\cdot)$ function. A piecewise polynomial interpolant provides dense output.

\subsection{Outline of paper}
We begin by introducing the geometry of the physical problem and some problem-specific background in \cref{sec_preamble}. The collocation scheme used is briefly described in \cref{sec_method}. The path-following method is described in \cref{sec_pf} for both the reduced and general problem using scattered data. In \cref{sec_accuracy}, we provide numerical results to determine the expected accuracy of the collocation and path-following methods. In \cref{sec_perf}, we evaluate the relative performance characteristics of the various methods and in \cref{sec_param} we describe how to choose optimal parameters. Finally, in \cref{sec_conclusions}, we provide some conclusions.

\section{Preamble}
\label{sec_preamble}

\subsection{Problem description}
Wave-current interaction problems are often studied by adopting a modal linear perturbative approach: waves are considered as first-order perturbations of a stationary, incompressible, and inviscid bulk fluid flow. In this context, waves are dispersive with the phase velocity of a wave dependent on the wave vector in a nonlinear manner. The relationship between the wave vector, $\vc{k}$, and the phase velocity, $\vc{c}$, is termed the \textit{dispersion relation} and is determined by factors such as water depth and background current.

In our context of first-order free-surface waves atop a vertical shear flow, the problem reduces to finding solutions of the eigenvalue problem formed from the Rayleigh instability equation and appropriate boundary conditions. The Rayleigh equation is a second order ODE that is equivalent to the Orr--Sommerfeld equation when viscous effects are neglected. A solution to the eigenproblem will yield a $\{k,c,w(z)\}$ triplet where $w(z)$ is the associated eigenfunction. The eigenfunction can be used to reconstruct the velocity and pressure field for the corresponding wave vector \cite[s. 4.4]{LiEllingsen2019}. Notably, there is substantial overlap with the literature on linear stability theory, e.g.\ \cite[ch.\ 2]{SchmidHenningson2001} or \cite{DrazinReid2004PLflow}, albeit with different boundary conditions.

Closed-form expressions for the dispersion relation for free-surface waves exist only in specific scenarios such as quiescent water \cite{MeiStiassnieYue2005} or a linear shear current (constant vorticity) \cite{Ellingsen2014a}. Integral approximation methods \cite{Skop1987}, \cite{KirbyChen1989}, \cite{EllingsenLi2017} and numerical schemes \cite{Shrira1993}, \cite{Zhang2005}, \cite{SmeltzerEllingsen2017}, \cite{LiEllingsen2019} exist to calculate the dispersion relation for arbitrary shear profiles.

\begin{figure}%
    \centering%
	\subfloat[An illustrative example of dispersive ring waves atop a shear current. An example shear profile is indicated with $z$-dependent arrows changing direction in the horizontal plane and an example $\vc{k}$ vector as $\vc{k}_0$.]{\includegraphics[trim=1.5cm 0.0cm 3.0cm 1.5cm,width=0.45\textwidth]{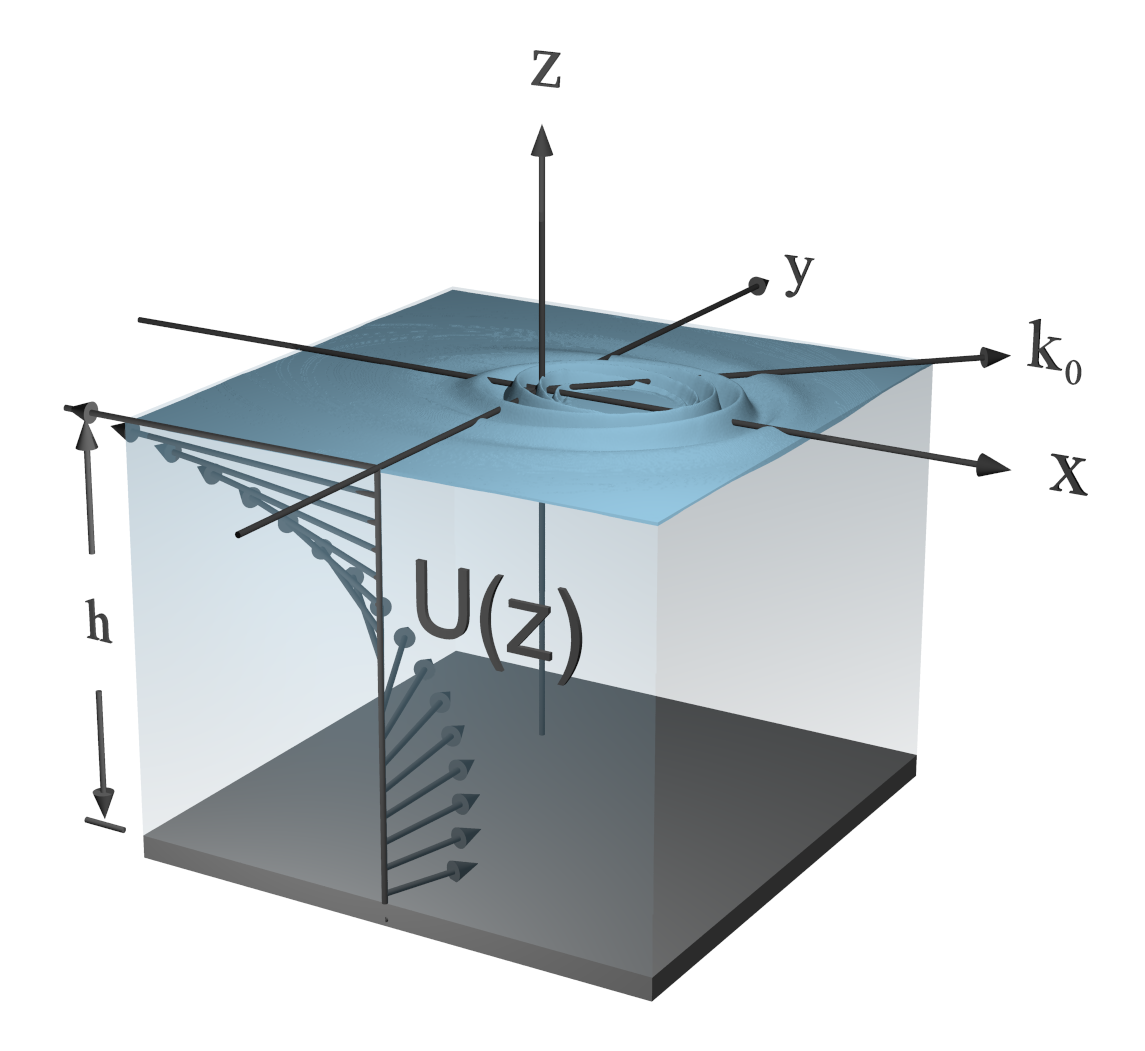} \label{fig_cp_diagram}}%
    \hfill
	\subfloat[Geometry used for reduced problem. Shear profile $U_\textrm{T}(z)$ shown in blue.]{\includegraphics[trim=1.5cm 0.0cm 1.0cm 1.0cm,width=0.45\textwidth]{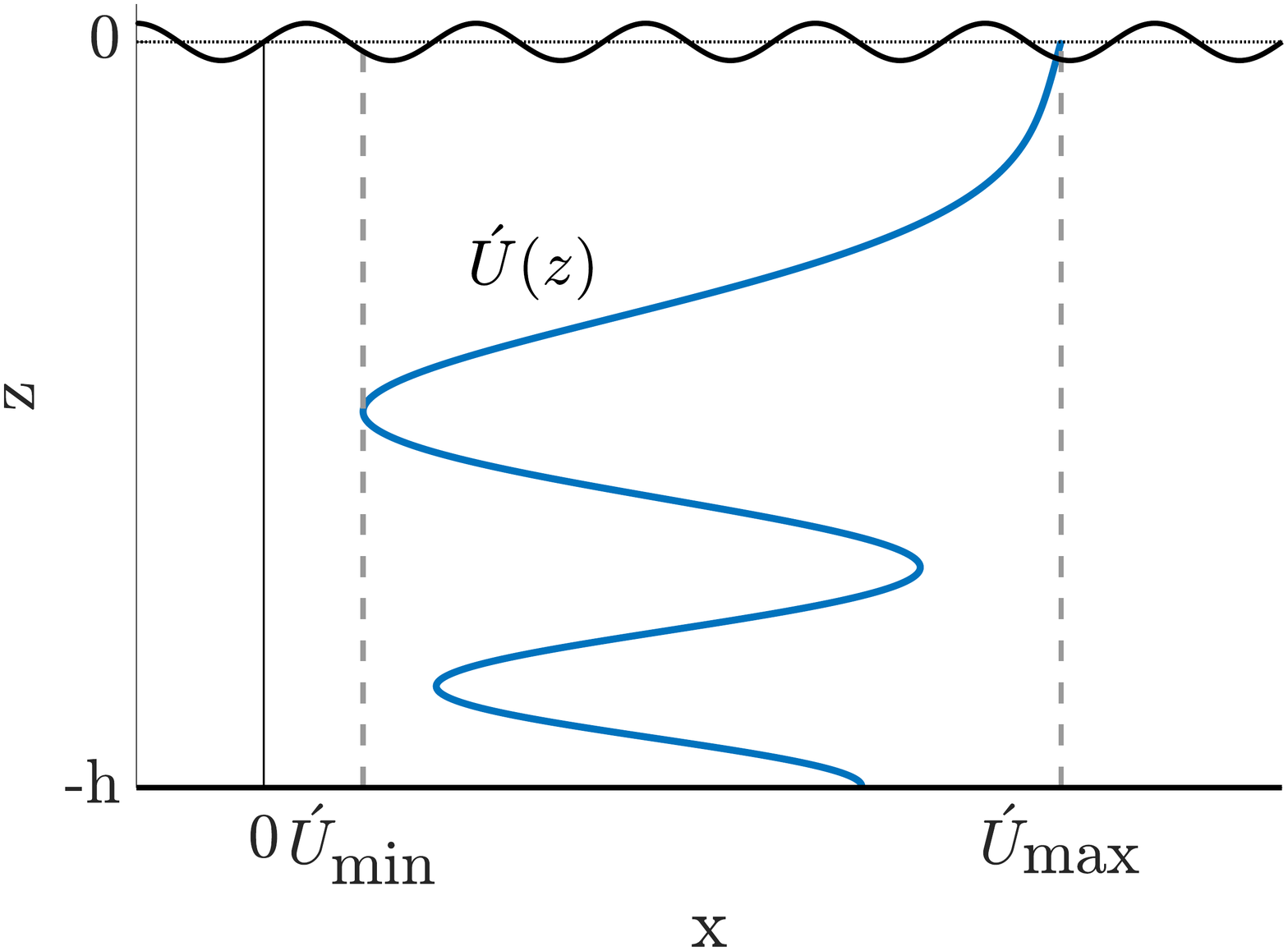} \label{fig_nondim_diagram}}%
    \caption{The problem setting.}      
    \label{fig_intro}
\end{figure}%

To make matters more concrete, we adopt the general approach used in \cite{Shrira1993}, \cite{LiEllingsen2019} and refer the reader to the derivations therein for full detail. For expediency, the approach is only summarised here. The physical model is depicted in \cref{fig_cp_diagram}.  Dimensional quantities are denoted with an acute, e.g.\ $\acute{k}$.

The background flow is specified by a shear profile $\acute{\vc{U}}(z) = ( \acute{U}_x(z), \acute{U}_y(z) )$: a 2-dimensional vector field describing the bulk fluid velocity in the horizontal plane for a given $\acute{z} \in [-\acute{h},0]$ where $\acute{h}$ is the constant fluid depth and the unperturbed surface is at $\acute{z}=0$. Let
\begin{equation}
\acute{U}_m := \max \left\{  \, | \sup_{z \in [0,-\acute{h}]} \acute{U}_x(z) - \inf_{z \in [0,-\acute{h}]} \acute{U}_x(z) |, \, | \sup_{z \in [0,-\acute{h}]} \acute{U}_y(z) - \inf_{z \in [0,-\acute{h}]} \acute{U}_y(z) | \, \right\}.
\end{equation}

We use $\acute{h}$ as a characteristic length scale and $\acute{U}_m$ as a characteristic velocity, to arrive at the following nondimensionalisation:
\begin{gather}
\vc{U}(z) = \acute{\vc{U}}(\acute{z}/\acute{h}) / \acute{U}_m, \quad z = \acute{z} / \acute{h},  \quad \vc{k} = \acute{\vc{k}} \acute{h}, \quad \vc{c} = \acute{\vc{c}} / \acute{U}_m, \quad h = 1,
\end{gather}
so that a notional shear strength can be expressed with Froude number, $F = \acute{U}_m / \sqrt{g \acute{h}}$. The shear profile must be suitably regular, so we impose that $U_x, U_y \in C^2([-1,0],\mathbb{R}) \cap L^2([-1,0],\mathbb{R})$. We also require $\vc{c}$ parallel to $\vc{k}$ and we henceforth only refer to scalar $c = |\vc{c}|$. It is assumed that the current can influence the waves but not conversely and, for clarity of exposition, we neglect surface tension. We adopt the Ansatz that perturbations are plane waves proportional to $\exp[i(\vc{k} \cdot \vc{x} - \omega t )]$ where the wave angular frequency is $\omega = kc$, and use a Fourier transform in the horizontal plane (coordinate space quantities indicated with a tilde),
\begin{equation}
\{ \tilde{u}, \tilde{v}, \tilde{w}, \tilde{p} \} = \frac{1}{(2\pi)^2} \int_{\mathbb{R}^2}  \{ u, v, w, p \} e^{i(\vc{k}\cdot\vc{x} - \omega t)} \, \drm \vc{k}.\label{eqn_ft}
\end{equation}
The velocity perturbations along the $x$, $y$, and $z$ axes are respectively $u=u(\vc{k},z)$, $v=v(\vc{k},z)$, and $w=w(\vc{k},z)$, whilst the pressure is $p=p(\vc{k},z)$. The governing equations are the linearised Euler equations and incompressibility condition, e.g.\ recall \cite{Shrira1993}:
\begin{align*}
ik_x u + ik_y v + w' &= 0, \quad
&i(\vc{k}\cdot\vc{U} - \omega) u + U_x'w &= -ik_xp/\rho, \\
i(\vc{k}\cdot\vc{U} - \omega) w &= -p'/\rho, \quad
&i(\vc{k}\cdot\vc{U} - \omega) v + U_y'w &= -ik_yp/\rho,
\end{align*}
with $k = |\vc{k}|$, and for shorthand $U_x=U_x(z)$, $U_y=U_y(z)$. After rearranging, the Rayleigh equation \cref{eqn_rayleigh3d_gov} is obtained which we write along with the relevant free-surface condition \cref{eqn_rayleigh3d_fs} (combined kinetic and dynamic boundary condition) and Dirichlet boundary condition for the fluid bottom \cref{eqn_rayleigh3d_bm},
\begin{subequations}
\label{eqn_rayleigh3d}
\begin{align}
w'' - \frac{\vc{k} \cdot \vc{U}''}{\vc{k} \cdot \vc{U} - kc} w - \mu w &= 0, \; &z \in (0,-h); \label{eqn_rayleigh3d_gov}\\
( \vc{k} \cdot \vc{U} - kc )^2 w' - [ ( \vc{k} \cdot \vc{U} - kc )\vc{k} \cdot \vc{U}' + F^{-2} k^2 ] w &= 0, \; &z = 0; \label{eqn_rayleigh3d_fs}\\
w &= 0, \; &z = -h; \label{eqn_rayleigh3d_bm}
\end{align}
\end{subequations}
for $\mu = k^2$. The velocity and pressure field for a specific wave vector can be recovered by substituting the appropriate eigenpair into \cite[s. 4.4]{LiEllingsen2019},
\begin{equation}
\begin{aligned}
(\vc{k}\cdot\vc{U} - \omega)w' -\vc{k}\cdot\vc{U}'w &= ik^2p/\rho,\\
ik_x[\vc{k}\cdot\vc{U}'w - (\vc{k}\cdot\vc{U} - \omega)w' ] - ik^2 U_x'w &= k^2(\omega - \vc{k}\cdot\vc{U})u,\\
ik_y[\vc{k}\cdot\vc{U}'w - (\vc{k}\cdot\vc{U} - \omega)w' ] - ik^2 U_y'w &= k^2(\omega - \vc{k}\cdot\vc{U})v.
\end{aligned}
\label{eqn_fullflow}
\end{equation}

The eigenvalue problem described by \cref{eqn_rayleigh3d} is inherently 1-dimensional. Usually, the physical problem is posed so that $U(z)$ is scalar for $\vc{U}(z) = (U(z),0)$, with scalar $k$ and $c$, see \cref{fig_nondim_diagram}.  However, there is  no difficulty in solving the more general physical problem by simply projecting $\vc{U}(z)$ along $\vc{k}$. Full 3-dimensional considerations only come to the fore when calculating the velocity and pressure field. To avoid the inherent ambiguities of `1d' / `2d' or `2d' / `3d' descriptions, we refer to the problem with scalar $U(z)$ as the \textit{reduced problem} and with vector $\vc{U}(z)$ as the \textit{general problem}. In \cref{subsec_gen_to_red}, we note that the reduced problem is equivalent to solving the general problem in a radial `slice' at some fixed angle $\theta_0$.

\begin{figure}%
    \centering%
    \subfloat[$\mu(c)$ spectrum for backward problem. Notice only one positive eigenvalue and series of negative eigenvalues distributed $\sim n^2$.]{{\includegraphics[trim=1.5cm 0.0cm 3.0cm 1.5cm,width=0.45\textwidth]{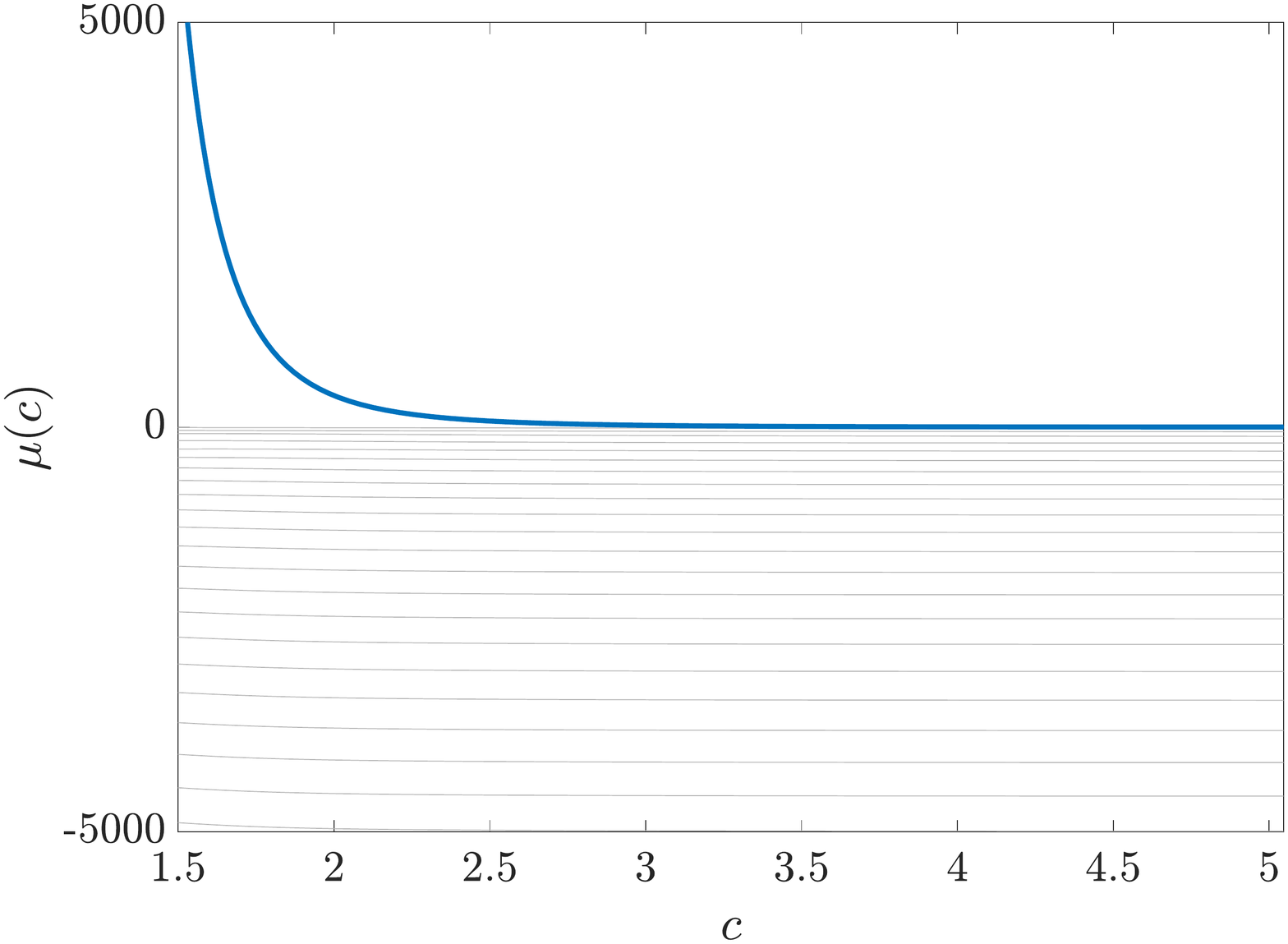}} \label{fig_spectrum_mu} }%
    \hfill
    \subfloat[Plot showing the sought positive and also negative solution branches along with the essential spectrum.]{{\includegraphics[trim=1.5cm 0.0cm 3.0cm 1.5cm,width=0.45\textwidth]{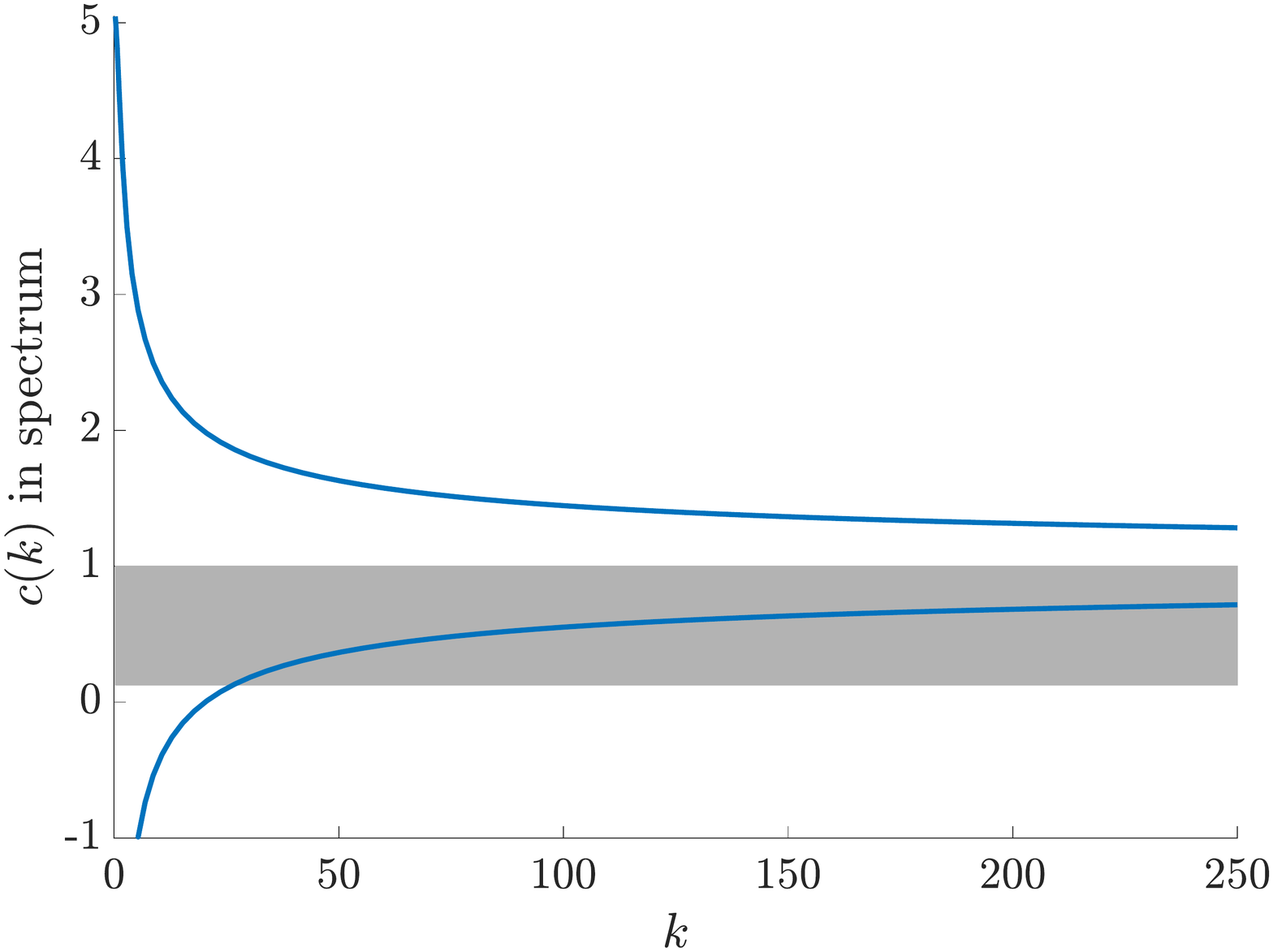} } \label{fig_spectrum_c}}%
    \caption{Parametrised plots of the real spectrum for both the backwards and forwards problem.}%
    \label{fig_spectrum_plots}%
\end{figure}

So far, we have deliberately avoided specifying which variable is the sought eigenvalue in \cref{eqn_rayleigh3d}: it can be chosen as either $\mu$ or $c$, with its counterpart parametrising the problem and always chosen to be real valued, in a similar manner to \cite[s.\ 7.4]{Boyd2001}. Since we are always choosing the parameter to be real-valued, we are concerned with a subset of the spectrum in each case and can plot this arrangement as a function of the parameter.
\begin{itemize}
\item The spectrum for $\mu(c)$, for $c$ in some suitable interval, is comprised of a countably infinite set of eigenvalues. The dominant eigenvalue, $\mu_1 = k^2$, in this case is the only positive eigenvalue, and corresponds to a propagating wave (for $\pm k$). The negative eigenvalues correspond to the countably infinite set of discrete $k$ arranged along the imaginary axis and are not mentioned further. See \cref{fig_spectrum_mu}.
\item The spectrum for $c(k)$, for $k > 0$, has both discrete and essential part ($c$ such that $U(z) - c = 0$, causing the ODE to become singular). In this case, the sought eigenvalue is again dominant but may be located within the same interval as the essential spectrum and therefore can be difficult to identify within numerical solution sets.  See \cref{fig_spectrum_c}.
\end{itemize}

\subsection{Problem types: forward, backward, and inverse}
\label{subsec_problem_types}
We distinguish three types of problem.
\begin{enumerate}
\item For shear profile, $U(z)$, and collection of wave numbers, $\{\sidx{k}{j}\}_{j=1}^J$, calculate associated velocities $\{c(\sidx{k}{j})\}_{j=1}^{J}$. We denote this as the \textit{forward problem}.
\item For shear profile, $U(z)$, and collection of velocities, $\{\sidx{c}{j}\}_{j=1}^J$, calculate associated wave numbers $\{k(\sidx{c}{j})\}_{j=1}^{J}$. We denote this as the \textit{backward problem}.
\item For collection of wave number and phase velocity pairs, $\{(\sidx{k}{j},\sidx{c}{j})\}_{j=1}^J$, determine shear profile, $U(z)$. We denote this as the \textit{inverse problem}, which is of an entirely different nature and herein not considered further.
\end{enumerate}

Both the forward and backward problem usually amount to calculating sufficiently many $\{ \sidx{k}{j},\sidx{c}{j} \}$ pairs as to adequately specify the full dispersion relation for a given shear profile.  For practical purposes, these problems are almost always posed as \cref{eqn_ft} with exponent of the form $i(kx - \omega(k) t)$ (see \cite[p. 77, eqn. 4.6]{Peregrine1976}, \cite[eqn. 1]{KirbyChen1989}, \cite[eqn. 2.4]{Shrira1993}, etc). This, by default, presupposes solving the forward problem.  Furthermore, handling of critical layers as in \cref{subsec_crit_c} is, from a numerical standpoint, easier for the forward problem.  Therefore, most of this paper concerns solution of the forward problem.

There are a few exceptions to this rule such as for wave problems with periodic or stationary time dependence, e.g.\ in ship waves. Hence, for purposes of completeness, we also describe solving the backwards problem using a rudimentary collocation scheme and the basic properties of the spectrum.

\subsection{Summary of numerical schemes}
For ease of reference, we denote the various numerical schemes used or described in this paper:
\begin{itemize}
\item \textit{CL-c} : Collocation scheme for the forward problem, see \cref{subsec_cl_sol_c}.
\item \textit{CL-k} : Collocation scheme for the backward problem, see \cref{subsec_cl_sol_k}.
\item \textit{PF-R-r-c} : A path-following scheme with dense output to solve the forward reduced problem along a fixed angle $\theta_0$ in the $\vc{k}$-plane, see \cref{subsubsec_pf_sys_c}.
\item \textit{PF-R-a-c} : A path-following scheme with dense output to solve the forward reduced problem for $c_{k_0}(\theta)$ along a fixed circle of radius $k = k_0$ with varying $\theta$ in the $\vc{k}$-plane, which we term the angular solution, see \cref{subsubsec_pf_sys_angular_c}.
\item \textit{PFmp-R-\{r,a\}-c} : An illustrative scheme using a single high precision QZ solve to improve accuracy of PF-R-\{r,a\}-c, detailed in \cref{subsubsec_improv_acc}.
\item \textit{PF-G-c} : A scheme which solves the forward general problem by using PF-R-r-c and PF-R-a-c to allow rapid interpolation with 2-dimensional scattered data query points in the $\vc{k}$-plane, see \cref{subsec_pf_general}.
\end{itemize}

The CL-c and CL-k schemes incur an eigenvalue calculation for each query point, so the computational cost will increase linearly with the number of query points. The arrangement of points in the $\vc{k}$-plane for the CL- schemes can be random without affecting computational cost.

The PF-R-r-c path-following algorithm is two stage: it first performs numerical integration to calculate control points along a radial `slice' at fixed angle $\theta_0$, which incurs a nominal initial computational cost; query points on that slice are then calculated using a Hermite interpolant. Although the computational cost of interpolation is linear in the number of query points, it is so light-weight as to be of almost negligible cost in most situations: so after the initial computation, very many query points can be calculated efficiently. The angular PF-R-a-c scheme is similar but instead calculates along a circular path at a fixed radius $k_0$.

The PF-G-c scheme is more involved because we accept query points in the $\vc{k}$-plane with no assumption on arrangement, i.e.\ scattered data. A naive approach would incur a complete first stage calculation of PF-R-r-c for every query point, which is unacceptable. The PF-G-c scheme instead precalculates a 2-dimensional polar grid of suitable control points and then can interpolate for query points at negligible cost.

Note that all methods presented can also make available the eigenfunction $w(z)$ so that the velocity and pressure field can be reconstructed using \cref{eqn_fullflow}.

\subsection{Existing algorithms and approximation methods}
\label{subsec_existing_alg}
Development in this area has been slightly unusual: despite the problem being readily amenable to numerical methods, there has been an emphasis on integral approximation schemes. In chronological order: Stewart \& Joy \cite{StewartJoy1974} (infinite depth), Skop \cite{Skop1987} (finite depth), Kirby \& Chen \cite{KirbyChen1989} (finite depth to 2nd order), and finally Ellingsen \& Li \cite{EllingsenLi2017}. Since our focus is on numerical methods, we do not address approximation schemes further.

The principal algorithm against which we compare is `DIM' from \cite{LiEllingsen2019}, which also contains a review of other numerical methods including the perennial piecewise-linear approach. For purposes of completeness, we also perform numerical simulations using a basic shooting method.

\subsection{Shear profiles and parameters used}
For later numerical tests, we define a test shear profile $U_\textrm{T}$ for the reduced problem, as shown in \cref{fig_nondim_diagram},
\begin{equation}
U_\textrm{T}(z) := \frac{\gamma}{2} ( 1 + \delta z ) \cos( \beta (-z)^\alpha ) + 1/2, \; \alpha = 2,\, \beta = 4\pi,\, \gamma = 1,\, \delta = 1/2.
\end{equation}
We choose the physical depth $\acute{h}=20$ and shear Froude number as $F^2 = 0.05$. We choose nondimensional $k \in I_k := [ \, 0.025, \, 250 \, ]$. This corresponds broadly to gravity waves in the air-water interface regime \cite[p.\ 4]{MeiStiassnieYue2005} with shortest period $\approx 0.2$s. The function is chosen as a suitable test candidate because it has several stationary points and cannot be approximated exactly over a finite dimensional polynomial basis.

For the figures produced from PF-G-c shown later in \cref{fig_interp_scatter}, we use $U_\textrm{T}(z)$ along the $x$ axis and an approximation of a flow in the Colombia River on the $y$ axis, which we denote $U_\textrm{CR}$. This is defined by a sixth degree polynomial and in our tests was scalled to have $F^2 \approx 0.01$; the precise definition of $U_\textrm{CR}$ is not so important as it is used only for the illustrative plots in \cref{subsec_pf_general}.

For more general choice of shear profile and parameters, it may be possible to create \textit{critical layers}. These are $z \in [-1,0]$ for which there exists some $c(k)$ such that $U(z) - c(k) = 0$, i.e\ depths $z$ for which the governing equation becomes singular. For our chosen shear profiles and parameters, critical layers are not encountered. Brief mention is made in \cref{subsec_crit_c} of how critical layers may be processed for CL-c.

\section{Collocation method for solving the dispersion relation}
\label{sec_method}

\subsection{General to reduced problem}
\label{subsec_gen_to_red}
The general problem \cref{eqn_rayleigh3d} can be simplified to a reduced problem by projecting $\vc{U}$ along $\vc{k}$, cf.\ \cite[p.\ 77]{Peregrine1976}\cite[p.\ 566]{Shrira1993}. Define the scalar shear profile for the reduced problem as $U_\theta(z) = (1/k) \vc{k}\cdot\vc{U} = \cos(\theta) U_x(z) + \sin(\theta) U_y(z)$ where $\theta$ is taken to be the standard angular polar coordinate for $\vc{k}$.

\subsection{Discretization of the equations}
We use sans serif notation to indicate matrices (uppercase) and vectors (lowercase), e.g. $\lamx{U}$ or $\lavc{w}$, to distinguish from their continuous counterparts.

Let $\elmidx{\zeta}{j} = \cos( (j-1) \pi / N_z)$, $j = 1 \ldots N_z+1$ be the Chebyshev--Gauss--Lobatto collocation points on $[-1,1]$ (second-kind points). We use the change of variable $z = (h/2)(\zeta - 1)$ to map $\elmidx{\zeta}{j}$ to $\elmidx{z}{j} \in [-1,0]$ in nondimensional coordinates and let $\lavc{z} := [ \sbns \elmidx{z}{1}, \ldots ,\elmidx{z}{N_z+1} \sbns ]^T$ be the associated column vector. Let $\ladiffmx$ be a corresponding square differentiation matrix (in practice, we calculate the first and second order differentiation matrices,  $\ladiffmx$ and $\ladiffsqmx$, using \texttt{poldif.m} from the Weideman--Reddy library \cite{WeidemanReddy2000} and then apply the `negative sum trick' as detailed in \cite{BaltenspergerTrummer2003}).  We define vector discretizations of the shear profile,
\begin{align*}
\lavc{u} = [ \elmidx{u}{1}, \ldots, \elmidx{u}{N_z+1} ]^T := [ U(\elmidx{z}{1}), \ldots, U(\elmidx{z}{N_z+1}) ]^T.
\end{align*}
Quantities $\lavc{u}'$ and $\lavc{u}''$ are similarly defined. Let $\lavc{w}:= w(\lavc{z})$. Define the diagonal matrices $\lamx{U} := \diag( \lavc{u} )$, $\lamx{U}' := \diag( \lavc{u}' )$, and $\lamx{U}'' := \diag( \lavc{u}'' )$.

The problem is a two-point boundary value problem so is amenable to the standard `row-replacement' strategy, see for example \cite[ch.\ 7]{Trefethen2000}. Specifically, we aim to construct eigenvalue equations which discretise the governing equation \cref{eqn_rayleigh3d_gov} using the `interior' rows 2 through $N_z$ of the differentiation and coefficient matrices. The free-surface boundary condition \cref{eqn_rayleigh3d_fs} is incorporated in the first row of the matrices. The bottom Dirichlet \cref{eqn_rayleigh3d_bm} boundary condition is accounted for by eliminating the last row and column in the matrices. For notational convenience, we define `interior' differentiation and shear profile matrices as $\ladiffmxint = \lamxidx{D}{lm}$, $\lamxint{U} = \lamxidx{U}{lm}$, $\lamxint{U}' = \lamxidx{U}{lm}'$, and $\lamxint{U}'' = \lamxidx{U}{lm}''$ with $l = 2 \ldots N_z$, $m = 1 \ldots N_z$ (in other words, eliminating the first and last rows, and last column). We also define a free-surface differentiation vector as the first row of $\ladiffmx$, $\lavcsb{d}{f} := \lamxidx{D}{1m}$, $m = 1 \ldots N_z$, again without the last column.

\subsection{Backward reduced problem (CL-k)}
\label{subsec_cl_sol_k}
Treating $c$ as a parameter and $k$ as the eigenvalue, we obtain a regular Sturm--Liouville problem on $z \in [-1,0]$ with eigenvalue $\mu = k^2$,
\begin{subequations}
\label{eqn_rayleigh2d_k_system}
\begin{align}
 \left( \frac{\drm^2}{\drm z^2} - \frac{U''(z)}{U(z)-c} \right) w(z) &=  \mu w(z), \; &z \in (-1,0); \label{eqn_rayleigh2d_k_gov}\\
(U-c)^2 w' - [ ( U - c )U' + F^{-2} ]w &= 0, \; &z = 0;  \label{eqn_rayleigh2d_k_fs}\\
w &= 0, \; &z=-1. \label{eqn_rayleigh2d_k_sb}
\end{align}
\end{subequations}

Let $\elmidx{q}{j}:= \elmidx{u}{j}'' / ( \elmidx{u}{j} - c )$, $\lavcint{q} := [ \elmidx{q}{2} \ldots \elmidx{q}{N_z} ]$, and $\lamxint{Q} = \diag(\lavcint{q})$. The discretization of \cref{eqn_rayleigh2d_k_gov} proceeds in the obvious manner,
\begin{equation}
\mu \lavc{w} = (  \ladiffsqmxint - \lamxint{Q} ) \lavc{w} = \lamx{R} \lavc{w}, \; \textrm{ for } \; \lamx{R}:= \ladiffsqmxint - \lamxint{Q}.
\end{equation}

Discretising \cref{eqn_rayleigh2d_k_fs} into a row vector gives,
\begin{equation}
\lavc{f}:= ( \elmidx{u}{0} - c )^2 \lavcsb{d}{f} - [ \sbns ( ( \elmidx{u}{0} - c ) \elmidx{u}{0}' + F^{-2} ), \sbns 0, \sbns \ldots, \sbns 0 \sbns ].
\end{equation}
Write
\begin{equation}
\lamx{A} = \begin{bmatrix}
\lavc{f} \\
\lamx{R}
\end{bmatrix}\textrm{, and } \lamx{B} = \diag( 0, 1, \ldots, 1 ),
\end{equation}
to obtain the generalised eigenvalue problem,
\begin{equation}
\label{eqn_eigs_ABk}
\lamx{A}\lavc{w} = \mu \lamx{B} \lavc{w}.
\end{equation}
Note that the only effect of $\lamx{B}$ is to ensure that the row of $\lamx{A}$ with the free-surface boundary condition is set equal to zero and is not dependent on the eigenvalue. \cref{eqn_eigs_ABk} can be solved in several ways, e.g.\ using MATLAB's implementation of QZ as \texttt{eig(A,B)}.

For a given $c$, there is a countably infinite set of discrete $\mu_j$ eigenvalues ordered $\mu_1 > \mu_2 > \ldots $. However, the only positive eigenvalue is $\mu_1$, which corresponds to the only real $k$, hence $\pm k$ represent the only propagating waves; we solve only for the positive branch. This is shown in \cref{fig_spectrum_mu}.

\subsection{Forward reduced problem (CL-c)}
\label{subsec_cl_sol_c}
Now, treating $k$ as a parameter and $c$ as the eigenvalue, we rewrite the reduced problem \cref{eqn_rayleigh2d_k_system} to emphasise the quadratic eigenvalue dependence in the free-surface boundary condition,
\begin{subequations}
\label{eqn_euler2d_c_system}
\begin{align}
\left( U \left( \frac{\drm^2}{\drm z^2} -k^2 \right) - U'' - c \left( \frac{\drm^2}{\drm z^2} -k^2 \right) \right) w &= 0, \; &z \in (-1,0); \label{eqn_euler2d_c_rayleigh}\\
c^2 w' +  c ( -2U w' + U' w ) + ( U^2 w' - U U' w - F^{-2} w ) &= 0, \; &z=0; \label{eqn_euler2d_c_fs}\\
w &= 0, \; &z=-1.\label{eqn_euler2d_c_sb}
\end{align}
\end{subequations}
We initially discretise \cref{eqn_euler2d_c_rayleigh} as
\begin{equation}
(  \lamxint{U} (\ladiffsqmxint - k^2 \lamx{I} ) - \lamxint{U}'' - c(\ladiffsqmxint - k^2 \lamx{I}) )\lavc{w} = 0.\label{eqn_dscr_c_v1}
\end{equation}
We proceed by expressing the free-surface condition as coefficients of the powers of $c$,
\begin{subequations}
\begin{align}
\lavcsb{f}{2} &= \lavcsb{d}{f}\\
\lavcsb{f}{1} &= -2 \elmidx{u}{0} \lavcsb{d}{f} + [ \sbns \elmidx{u}{0}', \sbns 0, \sbns \ldots, \sbns 0 \sbns ] \\
\lavcsb{f}{0} &= \elmidx{u}{0}^2 \lavcsb{d}{f} - [ \sbns \elmidx{u}{0}' \elmidx{u}{0} + F^{-2}, \sbns 0, \sbns \ldots, \sbns 0 \sbns ].
\end{align}
\end{subequations}
In the same manner, we now separate \cref{eqn_dscr_c_v1} into powers of $c$,
\begin{subequations}
\begin{align}
\lamxsb{R}{2} &= \lamx{0}\\
\lamxsb{R}{1} &= -( \ladiffsqmxint - k^2 \lamx{I} )\\
\lamxsb{R}{0} &= \lamxint{U} ( \ladiffsqmxint - k^2 \lamx{I} ) - \lamxint{U}''.
\end{align}
\end{subequations}
Define
\begin{equation}
\lamxsb{A}{2} = \begin{bmatrix}
\lavcsb{f}{2} \\
\lamxsb{R}{2}
\end{bmatrix}\textrm{, }
\lamxsb{A}{1} = \begin{bmatrix}
\lavcsb{f}{1} \\
\lamxsb{R}{1}
\end{bmatrix}\textrm{, and }
\lamxsb{A}{0} = \begin{bmatrix}
\lavcsb{f}{0} \\
\lamxsb{R}{0}
\end{bmatrix}.
\end{equation}
To obtain the sought solution, we solve the quadratic eigenproblem,
\begin{equation}
( c^2 \lamxsb{A}{2} + c \lamxsb{A}{1} + \lamxsb{A}{0} ) \lavc{w} = 0.\label{eqn_quad_eig_c}
\end{equation}

There are several techniques to solve the quadratic eigenvalue problem, although a direct linearisation and then using a QZ decomposition is sufficient in this setting. MATLAB's \texttt{polyeig(A2,A1,A0)} implements such a linearisation, although some care must be taken. In particular, the $\lamxsb{A}{2}$ matrix is badly rank-deficient. As a consequence, the QZ algorithm will return infinite and large-but-finite eigenvalues which are merely artefacts of the numerical method and must be removed.

The spectrum has two discrete branches and essential spectrum; we seek the positive branch (greatest eigenvalue), which corresponds to propagating waves. The essential spectrum contains all $c$ such that $U(z) - c = 0$ for some $z \in [-1,0]$, see \cite{DrazinReid2004PLflow} and \cref{fig_spectrum_c}.

\subsection{Solving with critical layers in the forward problem}
\label{subsec_crit_c}
Let
\begin{equation*}
U_{\min} := \inf_{z \in [-1,0]} \{ U(z) \} \textrm{ and } U_{\max} := \sup_{z \in [-1,0]} \{ U(z) \}.
\end{equation*}
A critical layer exists if phase velocity $c \in [U_{\min}, U_{\max}]$ (see shaded region in \cref{fig_spectrum_c}), in other words if $c$ is in the region occupied by the essential spectrum. The QZ algorithm returns many points from the essential spectrum and, if $c \in [U_{\min}, U_{\max}]$, they are numerically indistinguishable. However, the eigenvectors from the essential spectrum have singular behaviour in the interior of their domain whereas this is not true for the eigenvector corresponding to sought eigenvalue $c$. Thus, the sought eigenvalue can, in principle at least, be identified and computation may still proceed when critical layer(s) are present. The question of critical layers is, however, not central to the theme of this paper and so shall not be mentioned further.

\section{Path-following method for calculating the dispersion relation curve}
\label{sec_pf}

\subsection{Review of Loisel--Maxwell path-following method for the field of values}
\label{subsec_lm_pf}
In \cite{LoiselMaxwell2018}, the authors describe a path-following method to calculate the field of values boundary of a matrix, which we now briefly summarise. It concerns the solution of a parametrised Hermitian eigenvalue problem (which bounds the projection of the field of values onto the real axis),
\begin{equation}
\fnherm{e^{i\tau} \lamx{A} }\lavc{u}(\tau) = \lambda(\tau) \lavc{u}(\tau) \textrm{ for } \lamx{A} \in \mtxsqrcplx{n},\, \lavc{u} \in \mathbb{C}^n,\, \lambda \in \mathbb{R},\, \tau \in [0,2\pi),\label{eqn_fov_herm}
\end{equation}
taking $(\lambda(\tau), \lavc{u}(\tau))$ as the dominant eigenpair where $\fnherm{\lamx{A}}:= (1/2)(\lamx{A} + \lamx{A}^{*})$ is the Hermitian part, $\fnsherm{\lamx{A}}:= (1/2)( \lamx{A} - \lamx{A}^{*})$ is the skew-Hermitian part of the given matrix $\lamx{A}$, and $\lamx{A}^{*}$ is the conjugate-transpose. Here, and in the remainder of the paper, the overdot notation is used to indicate derivatives with respect to the problem parameter. This is to emphasise the parameter-varying or ``time-varying'' nature of the problems.

Note that \cref{eqn_fov_herm} is well-defined except perhaps for a finite number of $\tau_j$ due to elementwise analyticity of the elements of $\fnherm{e^{i\tau} \lamx{A} }$ and the analyticity, up to ordering, of the eigenvalue and eigenvectors. Differentiating \cref{eqn_fov_herm} gives,
\begin{equation}
\fnherm{e^{it} \lamx{A}}\lavcd{u}(\tau) - \dot{\lambda}(\tau)\lavc{u}(\tau) - \lambda(\tau)\lavcd{u}(\tau) = -i\fnsherm{e^{i\tau} \lamx{A} }\lavc{u}(\tau).
\end{equation}
Since the system is under-determined, an additional constraint that $\lavc{u}(\tau)$ must be tangent to its (elementwise) derivative is included, giving the system,
\begin{equation}
\begin{cases}
\fnherm{e^{it} \lamx{A}}\lavcd{u}(\tau) - \dot{\lambda}(\tau)\lavc{u}(\tau) - \lambda(\tau)\lavcd{u}(\tau) = -i\fnsherm{e^{it} \lamx{A} }\lavc{u}(\tau)\\
\lavc{u}(\tau)^*\lavcd{u}(\tau) = 0,
\end{cases}
\end{equation}
which can be rewritten in matrix form,
\begin{equation}
\label{eqn_fov_system}
\mtrxtwobytwo{\fnherm{e^{i\tau} \lamx{A} } - \lambda(\tau)I}{-\lavc{u}(\tau)}{-\lavc{u}(\tau)^*}{0} \mtrxtwobyone{\lavcd{u}(\tau)}{\dot{\lambda}(\tau)} = \mtrxtwobyone{-i\fnsherm{e^{i\tau} \lamx{A} }\lavc{u}(\tau)}{0}.
\end{equation}
The system described by \cref{eqn_fov_system} can be solved for $[ \sbns \lavcd{u}(\tau)^{*} \sbns \dot{\lambda}(\tau)^{*} \sbns ]^{*}$ using a linear solver and used as the $F(\cdot)$ function for a Runge--Kutta numerical integrator, which then generates control points along the curve. The authors use the Dormand--Prince RK5(4)7M method \cite[p.\ 23]{DormandPrince1980} and interpolation method of Shampine \cite[p.\ 148]{Shampine1986}. The near-interpolant solution from this method is $5^\textrm{th}$ order accurate.

\subsection{Path-following method for forward reduced problem}
\label{subsec_pf_reduced}
We now extend the same process to the quadratic eigenvalue problem posed in \cref{subsec_cl_sol_c}. Recall \cref{eqn_quad_eig_c},
\begin{equation}
\Bigl( c^2(k) \lamxsb{A}{2}(k) + c(k) \lamxsb{A}{1}(k) + \lamxsb{A}{0}(k) \Bigr) \lavc{w}(k) = 0,
\end{equation}
which upon differentiating (indicated with overdot) with respect to $k$ gives,
\begin{multline}
\Bigl( c^2(k) \lamxsbd{A}{2}(k) + c(k) \lamxsbd{A}{1}(k) + \lamxsbd{A}{0}(k) +  2c(k)\dot{c}(k) \lamxsb{A}{2}(k) + \dot{c}(k) \lamxsb{A}{1}(k) \Bigr) \lavc{w}(k)+ \\
\Bigl( c^2(k) \lamxsb{A}{2}(k) + c(k) \lamxsb{A}{1}(k) + \lamxsb{A}{0}(k) \Bigr) \lavcd{w}(k) = 0.
\end{multline}
We further impose that $\lavc{w}(k)^{*} \lavcd{w}(k) = 0$. Writing in matrix form,
\begin{multline}
\begin{bmatrix}
\Bigl(  c^2(k) \lamxsb{A}{2}(k) + c(k) \lamxsb{A}{1}(k) + \lamxsb{A}{0}(k) \Bigr) & \Bigl( 2c(k) \lamxsb{A}{2}(k) + \lamxsb{A}{1}(k) \Bigr) \lavc{w}(k) \\
\lavc{w}^{*}(k) & 0
\end{bmatrix}
\begin{bmatrix}
\lavcd{w}(k) \\
\dot{c}(k)
\end{bmatrix}\\
 = \begin{bmatrix}
-\Bigl( c^2(k) \lamxsbd{A}{2}(k) + c(k) \lamxsbd{A}{1}(k) + \lamxsbd{A}{0}(k) \Bigr) \lavc{w}(k) \\ 0
\end{bmatrix}.\label{eqn_pf_general}
\end{multline}
This is the general form in which the structure is clear. In the subsections below, we perform the same derivation but include the specific expressions for the radial and angular paths including boundary conditions.

The approach taken is analogous to \cite{LoiselMaxwell2018}: an initial eigenpair $\{ \elminitval{c}, \lavcinitval{w} \}$ is calculated using CL-c. Then by using \cref{eqn_pf_general} to solve for $[ \sbns \lavcd{w}(k)^{*} \sbns \dot{c}(k)^{*} \sbns ]^{*}$, numerical integration can proceed along the curve in both directions. Hermite interpolation can then be used to query at arbitrary $k$.

\subsubsection{System of equations along radial slice at fixed $\theta$ (PF-R-r-c)}
\label{subsubsec_pf_sys_c}
For PF-R-r-c, we fix angle $\theta = \elminitval{\theta}$ and parametrise by $k$. Thus, we are in the setting of the reduced problem with the constant shear profile being the relevant reduced shear profile, $U_\theta(z)$. Writing \cref{eqn_rayleigh2d_k_gov} in matrix form with $c$ as eigenvalue and explicit dependence on parameter $k$,
\begin{equation}
\Bigl( \lamxint{U} - c(k) \lamxint{I} \Bigr) \ladiffsqmxint \lavc{w}(k) - \lamxint{U}'' \lavc{w}(k) = k^2 \Bigl( \lamxint{U} - c(k) \lamxint{I} \Bigr) \lavc{w}(k). \label{eqn_rayleigh2d_PF_c_gov}
\end{equation}
For notational succinctness, we use the shorthand $c=c(k)$ and $\lavc{w}=\lavc{w}(k)$. Differentiating \cref{eqn_rayleigh2d_PF_c_gov} with respect to $k$ (indicated by an overdot) gives,
\begin{equation}
( \lamxint{U} - c \lamxint{I} ) \ladiffsqmxint \lavcd{w} - \dot{c} \ladiffsqmxint \lavc{w}  - \lamxint{U}'' \lavcd{w} =  2k ( \lamxint{U}  - c \lamxint{I} ) \lavc{w}  - \dot{c} k^2 \lavc{w} + k^2 ( \lamxint{U} - c \lamxint{I} ) \lavcd{w}
\end{equation}

The free-surface condition can be written as,
\begin{equation}
( \elmidx{u}{0} - c )^2 \lavcsb{d}{f} \lavc{w}
+ [ \sbns ( -\elmidx{u}{0} \elmidx{u}{0}' + c \elmidx{u}{0}' - F^{-2} ), \sbns 0, \sbns \ldots, \sbns 0 \sbns ]\lavc{w} = 0.\label{eqn_rayleigh2d_PF_c_fs}
\end{equation}
Differentiating \cref{eqn_rayleigh2d_PF_c_fs} with respect to $k$,
\begin{multline}
-2 \dot{c} ( \elmidx{u}{0} - c ) \lavcsb{d}{f} \lavc{w}
+ ( \elmidx{u}{0} - c )^2  \lavcsb{d}{f} \lavcd{w}\\
+ [ \sbns ( \dot{c} \elmidx{u}{0}' ), \sbns 0, \sbns \ldots, \sbns 0 \sbns ]\lavc{w} 
+ [ \sbns -( \elmidx{u}{0}' ( \elmidx{u}{0} - c ) + F^{-2} ), \sbns 0, \sbns \ldots, \sbns 0 \sbns ] \lavcd{w} = 0.
\end{multline}

Upon rearranging terms, we define the block matrices:
\begin{subequations}
\begin{gather}
\lamx{P}(k,c) := \begin{bmatrix}
( \elmidx{u}{0} - c )^2 \lavcsb{d}{f} + [ \sbns -( \elmidx{u}{0}' ( \elmidx{u}{0} - c ) + F^{-2} ), \sbns 0, \sbns \ldots, \sbns 0 \sbns ] \\
( \lamxint{U} - c \lamxint{I} ) \ladiffsqmxint  - \lamxint{U}'' - k^2 ( \lamxint{U} - c \lamxint{I} ) 
\end{bmatrix},\\
\lamx{Q}(k,c) := \begin{bmatrix}
-2 ( \elmidx{u}{0} - c ) \lavcsb{d}{f} + [ \sbns ( \elmidx{u}{0}' ), \sbns 0, \sbns \ldots, \sbns 0 \sbns ] \\
-\ladiffsqmxint + k^2 \lamxint{I}
\end{bmatrix}, \textrm{and}\\
\lamx{R}(k,c) := \begin{bmatrix}
0 \\
 2 k ( \lamxint{U} - c \lamxint{I} ) 
\end{bmatrix},
\end{gather}\label{eqn_p_q_r}
\end{subequations}
so that the system of ODEs can now be written in matrix form as:
\begin{equation}
\label{eqn_PF_c_system}
\underbrace{\begin{bmatrix}
\lamx{P}(k,c) & \lamx{Q}(k,c) \lavc{w}(k) \\
-\lavc{w}^*(k) & 0
\end{bmatrix}}_{\lamx{M}(k,c,\lavc{w})}
\underbrace{\begin{bmatrix}
\lavcd{w}(k) \\
\dot{c}(k)
\end{bmatrix}}_{\lavc{v}(k)}
= \underbrace{\begin{bmatrix}
\lamx{R}(k,c)\lavc{w}(k) \\
0
\end{bmatrix}}_{\lavc{b}(k,c,\lavc{w})}.
\end{equation}
Note that we do not include the row corresponding to the bottom surface, only the free-surface is included.

\subsubsection{System of equations along angular circle at fixed $k$ (PF-R-a-c)}
\label{subsubsec_pf_sys_angular_c}
For PF-R-r-c, the angle $\theta$ and hence the shear profile was held constant. For PF-R-a-c, we instead hold $k$ constant and seek to use a $\theta$ angular dependence. Therefore, we must also specify the parametrisation of the shear profile.
\begin{equation*}
U_k(\theta,z) = \cos(\theta) U_x(z) + \sin(\theta) U_y(z).
\end{equation*}
So that in matrix form,
\begin{equation*}
\lamxidx{U}{}(\theta) = \cos(\theta)\, \lamxux + \sin(\theta)\, \lamxuy,
\end{equation*}
and, upon differentiation with respect to $\theta$ (indicated with an overdot),
\begin{equation*}
\lamxdf{U}{}{}(\theta) = -\sin(\theta)\, \lamxux + \cos(\theta)\, \lamxuy.
\end{equation*}

Our starting point is the same, we use \cref{eqn_rayleigh2d_PF_c_gov} but instead hold $k$ constant and take the derivative with respect to $\theta$. Temporarily adopting the abbreviated notation $\lamxint{U} = \lamxint{U}(\theta)$, $\lavc{w} = \lavc{w}(\theta)$, and $c = c(\theta)$: 
\begin{multline}
( \lamxintd{U} - \dot{c} \lamxint{I} ) \ladiffsqmxint \lavc{w} + ( \lamxint{U} - c \lamxint{I} ) \ladiffsqmxint \lavcd{w} - \lamxintd{U}'' \lavc{w} - \lamxint{U}'' \lavcd{w} \\
 = k^2 ( ( \lamxintd{U}  - \dot{c} \lamxint{I} ) \lavc{w} + ( \lamxint{U} - c \lamxint{I} ) \lavcd{w} ),
\end{multline}
As before, the free-surface condition is \cref{eqn_rayleigh2d_PF_c_fs}, which we take the derivative of with respect to $\theta$ using the shorthand $\elmidxd{u}{0} = \elmidxd{u}{0}(\theta)$,
\begin{multline}
2( \elmidx{u}{0} - c ) ( \elmidxd{u}{0} - \dot{c} ) \lavcsb{d}{f} \lavc{w} + ( \elmidx{u}{0} - c )^2 \lavcsb{d}{f} \lavcd{w} \\
+ [ \sbns -(  \elmidxd{u}{0}' ( \elmidx{u}{0} - c ) + \elmidx{u}{0}' ( \elmidxd{u}{0} - \dot{c} ) ), \sbns 0, \sbns \ldots, \sbns 0 \sbns ]\lavc{w}  \\
+ [ \sbns -( \elmidx{u}{0}' ( \elmidx{u}{0} - c ) + F^{-2} ), \sbns 0, \sbns \ldots, \sbns 0 \sbns ] \lavcd{w} = 0.
\end{multline}

In a similar manner to before, we define the block matrices:
\begin{subequations}
\begin{gather}
\lamx{P}(\theta,c) := \begin{bmatrix}
( \elmidx{u}{0} - c )^2 \lavcsb{d}{f}  + [ \sbns -( \elmidx{u}{0}' ( \elmidx{u}{0} - c ) + F^{-2} ), \sbns 0, \sbns \ldots, \sbns 0 \sbns ] \\
( \lamxint{U} - c \lamxint{I} ) \ladiffsqmxint - \lamxint{U}'' - k^2 ( \lamxint{U} - c \lamxint{I} )
\end{bmatrix},\\
\lamx{Q}(\theta,c) := \begin{bmatrix}
-2 ( \elmidx{u}{0} - c ) \lavcsb{d}{f} + [ \sbns ( \elmidx{u}{0}' ), \sbns 0, \sbns \ldots, \sbns 0 \sbns ] \\
-\ladiffsqmxint + k^2 \lamxint{I}
\end{bmatrix}, \textrm{and}\\
\lamx{R}(\theta,c) := \begin{bmatrix}
- 2 ( \elmidx{u}{0} - c ) \elmidxd{u}{0} \lavcsb{d}{f} + [ \sbns ( \elmidxd{u}{0} \elmidx{u}{0}' + \elmidxd{u}{0}' ( \elmidx{u}{0} - c ) ), \sbns 0, \sbns \ldots, \sbns 0 \sbns ]  \\
-\lamxintd{U} \ladiffsqmxint + \lamxintd{U}'' + k^2 \lamxintd{U}
\end{bmatrix},
\end{gather}
\end{subequations}
so that the system of ODEs can now be written in matrix form as:
\begin{equation}
\label{eqn_PF_angular_c_system}
\underbrace{\begin{bmatrix}
\lamx{P}(\theta,c) & \lamx{Q}(\theta,c) \lavc{w}(\theta) \\
-\lavc{w}^*(\theta) & 0
\end{bmatrix}}_{\lamx{M}(\theta,c,\lavc{w})}
\underbrace{\begin{bmatrix}
\lavcd{w}(\theta) \\
\dot{c}(\theta)
\end{bmatrix}}_{\lavc{v}(\theta)}
= \underbrace{\begin{bmatrix}
\lamx{R}(\theta,c)\lavc{w}(\theta) \\
0
\end{bmatrix}}_{\lavc{b}(\theta,c,\lavc{w})}.
\end{equation}

Note that the $\lamx{P}$ and $\lamx{Q}$ matrix have the same structure as in \cref{eqn_PF_c_system}, it is $\lamx{R}$ that changes.

\subsubsection{Path-following algorithm specification for reduced problem}
We describe the algorithm for PF-R-r, the algorithm for PF-R-a follows in the obvious manner. Using the definitions of $\lamx{P}$, $\lamx{Q}$, $\lamx{R}$ from \cref{eqn_p_q_r} define matrix and vector functions,
\begin{equation}
\lamx{M}(k,c,\lavc{w}) := \begin{bmatrix}
\lamx{P}(k,c) & \lamx{Q}(k,c) \lavc{w}(k) \\
-\lavc{w}^*(k) & 0
\end{bmatrix}, \qquad \lavc{b}(k,c,\lavc{w}) := \begin{bmatrix}
\lamx{R}(k,c)\lavc{w}(k) \\
0
\end{bmatrix}.
\end{equation}
Given a candidate $\lavc{v}(k):= [ \sbns \lavc{w}(k)^{*} \sbns c(k)^{*} \sbns ]^{*}$, define the Runge--Kutta $F(\cdot)$ function as,
\begin{equation}
F\left(k,\begin{bmatrix}
\lavc{w}(k) \\
c(k)
\end{bmatrix}\right) = \begin{bmatrix}
\lavcd{w}(k) \\
\dot{c}(k)
\end{bmatrix} = \lamx{M}(k,c,\lavc{w})^{-1} \, \lavc{b}(k,c,\lavc{w}).
\end{equation}
$F(\cdot)$ can be easily obtained using a linear solver, such as LU decomposition.

\begin{figure}%
    \centering%
    \subfloat[Dormand--Prince control points indicated by blue circles.]{{\includegraphics[trim=1.5cm 0.0cm 3.0cm 1.5cm,width=0.45\textwidth]{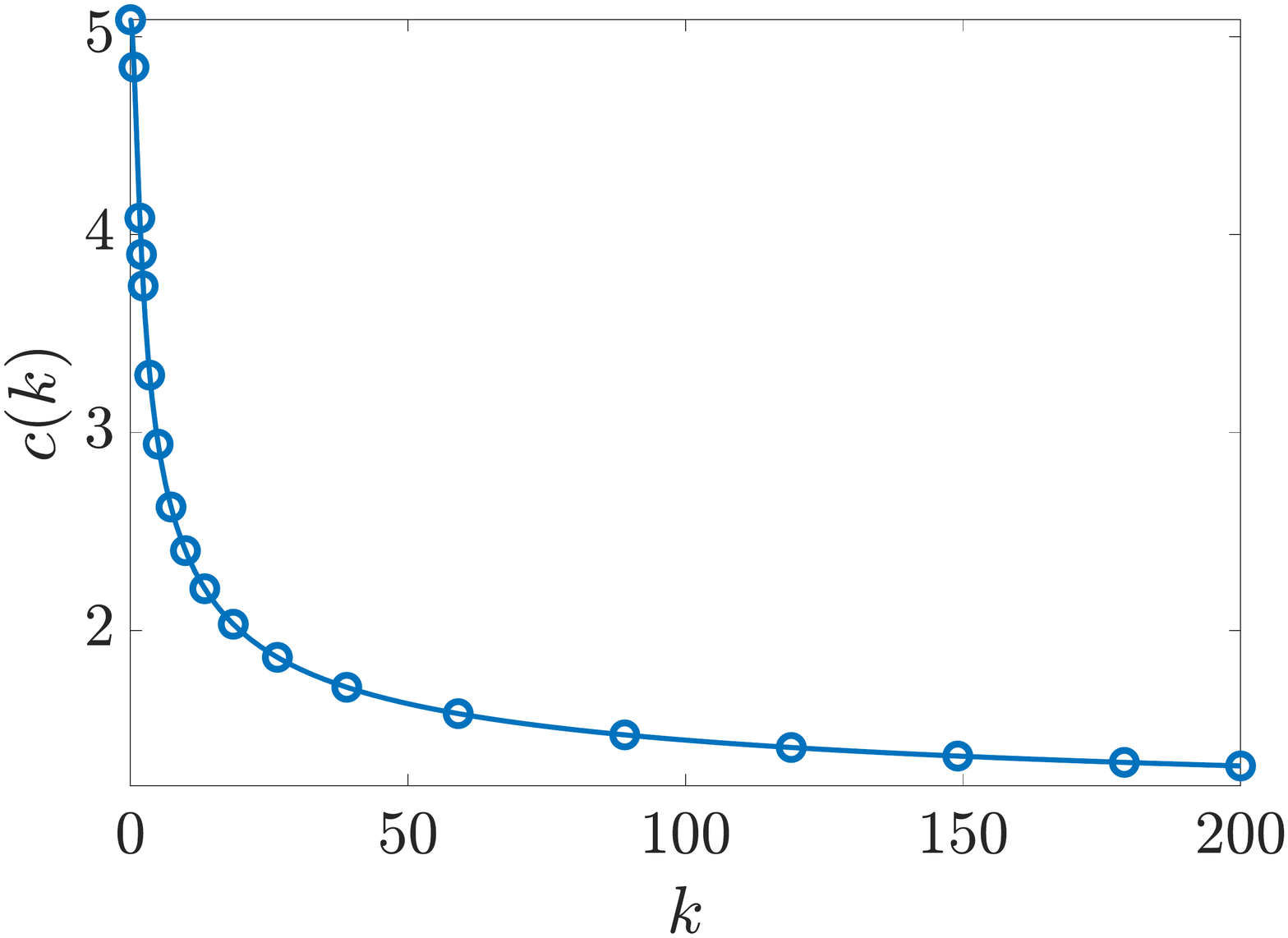}} \label{fig_interp1d_normal} }%
    \hfill
    \subfloat[Zoomed-in. Sample interpolant query points shown with red asterisks.]{{\includegraphics[trim=1.5cm 0.0cm 3.0cm 1.5cm,width=0.45\textwidth]{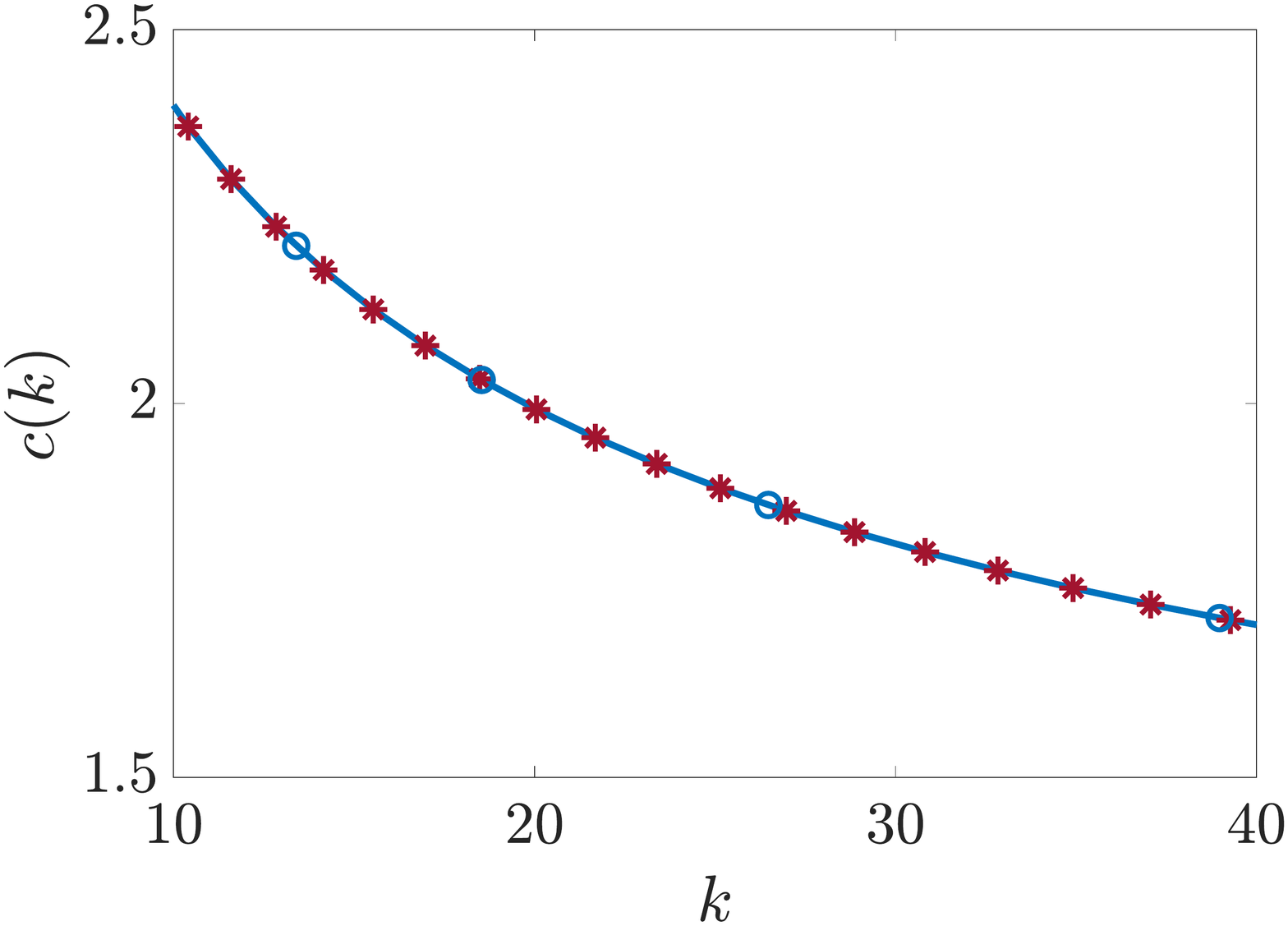} } \label{fig_interp1d_zoom}}%
    \caption{Numerical integration of dispersion relation curve for shear profile $U_\textrm{T}$. Tolerance for integrator was $10^{-6}$.}%
    \label{fig_interp1d}%
\end{figure}

The algorithm requires an initial $\lavcinitval{v} = \lavc{v}(\elminitval{k})$ calculated using CL-c. As in \cite{LoiselMaxwell2018}, the Dormand--Prince RK5(4)7M method \cite[p. 23]{DormandPrince1980} and Hermite interpolation strategy of Shampine \cite[p. 148]{Shampine1986} is used. We use automatic stepsize control as described in \cite[p. 167]{HairerNorsettWanner1993}. For an interval $[ \sidx{k}{j}, \sidx{k}{j+1} ]$ with midpoint $\sidx{k}{\textrm{mid}}$, the integrator produces control points $\{ \sidx{\lavc{v}}{j}, \sidx{\lavc{v}}{j}, \sidx{\lavc{v}}{j+1/2}, \sidx{\lavc{v}}{j+1}, \sidx{\lavc{v}}{j+1} \}$ where $\sidx{\lavc{v}}{j} = \lavc{v}(\sidx{k}{j})$ and $\sidx{\lavc{v}}{j+1/2} = \lavc{v}(\sidx{k}{\textrm{mid}})$. Thus, after numerical integration, a solution set of $\sidx{\lavc{v}}{j}$, $\sidx{\lavc{v}}{j}$, and $\sidx{\lavc{v}}{j+1/2}$ is obtained upon which piecewise Hermite interpolation can be performed. If both $c(k)$ and the eigenvector $\lavc{w}(k)$ is required then interpolation is over $N+1$ length vectors; if only $c(k)$ is required then interpolation is only 1-dimensional. Example output is shown in \cref{fig_interp1d} (1-dimensional output).

\subsubsection{PFmp-R-\{r,a\}-c: improving accuracy for PF-R-\{r,a\}-c}
\label{subsubsec_improv_acc}
As shall be described in \cref{sec_accuracy}, the error in the CL- methods are determined almost entirely by roundoff error incurred during the solution of the quadratic eigenvalue problem in \texttt{double} precision. The path-following algorithm essentially maintains the same error as is present in the initial $\lavcinitval{v}$.  By calculating $\lavc{v}(\elminitval{k})$ in high precision arithmetic then executing the path-following schemes in \texttt{double} precision as normal, an improvement in accuracy of two to three orders of magnitude is obtained. This is discussed further in \cref{subsec_accuracy_eigenvalues_nz}.

\subsection{Path-following method for forward general problem (PF-G-c)}
\label{subsec_pf_general}
The PF-R-r and PF-R-a algorithms can be combined to create an efficient algorithm that can process scattered data query points, which we describe below.

\begin{figure}%
    \centering%
    \subfloat[Planar plot of PF-R-a used at normial radius $\elminitval{k}$, interpolated at angles $\sidx{\theta}{j}$. Blue circles are control points, red astrisks are interpolation points. Angles indicated in dotted grey.]{{\includegraphics[trim=1.5cm 0.0cm 3.0cm 1.5cm,width=0.45\textwidth]{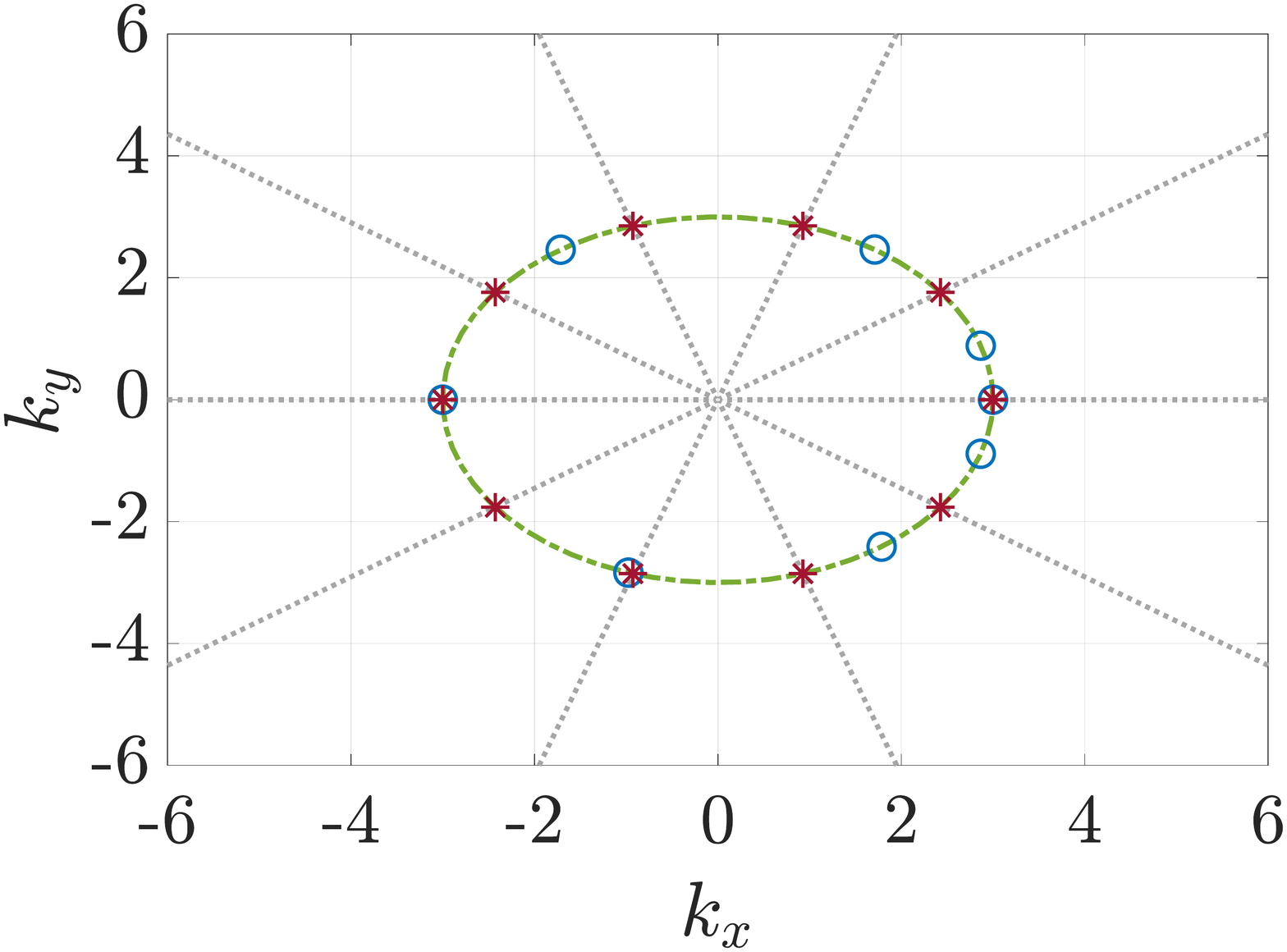}} \label{fig_interp_precomp_ang_a} }%
    \hfill
    \subfloat[3d plot as panel (a).]{{\includegraphics[trim=1.5cm 0.0cm 3.0cm 1.5cm,width=0.45\textwidth]{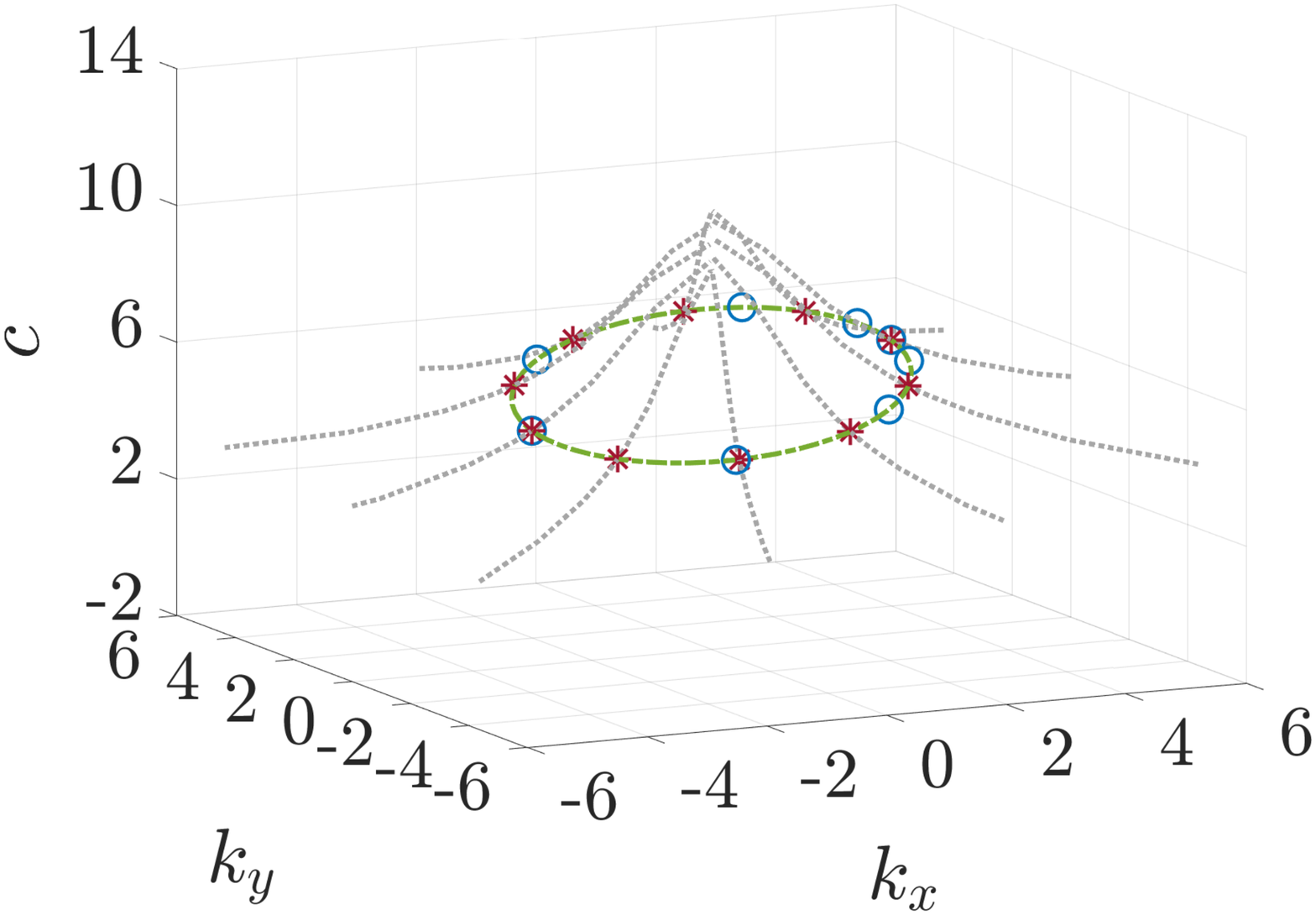} } \label{fig_interp_precomp_ang_b}}\\%
    \subfloat[Planar plot of PF-R-r used along each $\sidx{\theta}{j}$. Blue circles are control points, red astrisks are interpolation points (the $\sidx{k}{i}$).]{{\includegraphics[trim=1.5cm 0.0cm 3.0cm 1.5cm,width=0.45\textwidth]{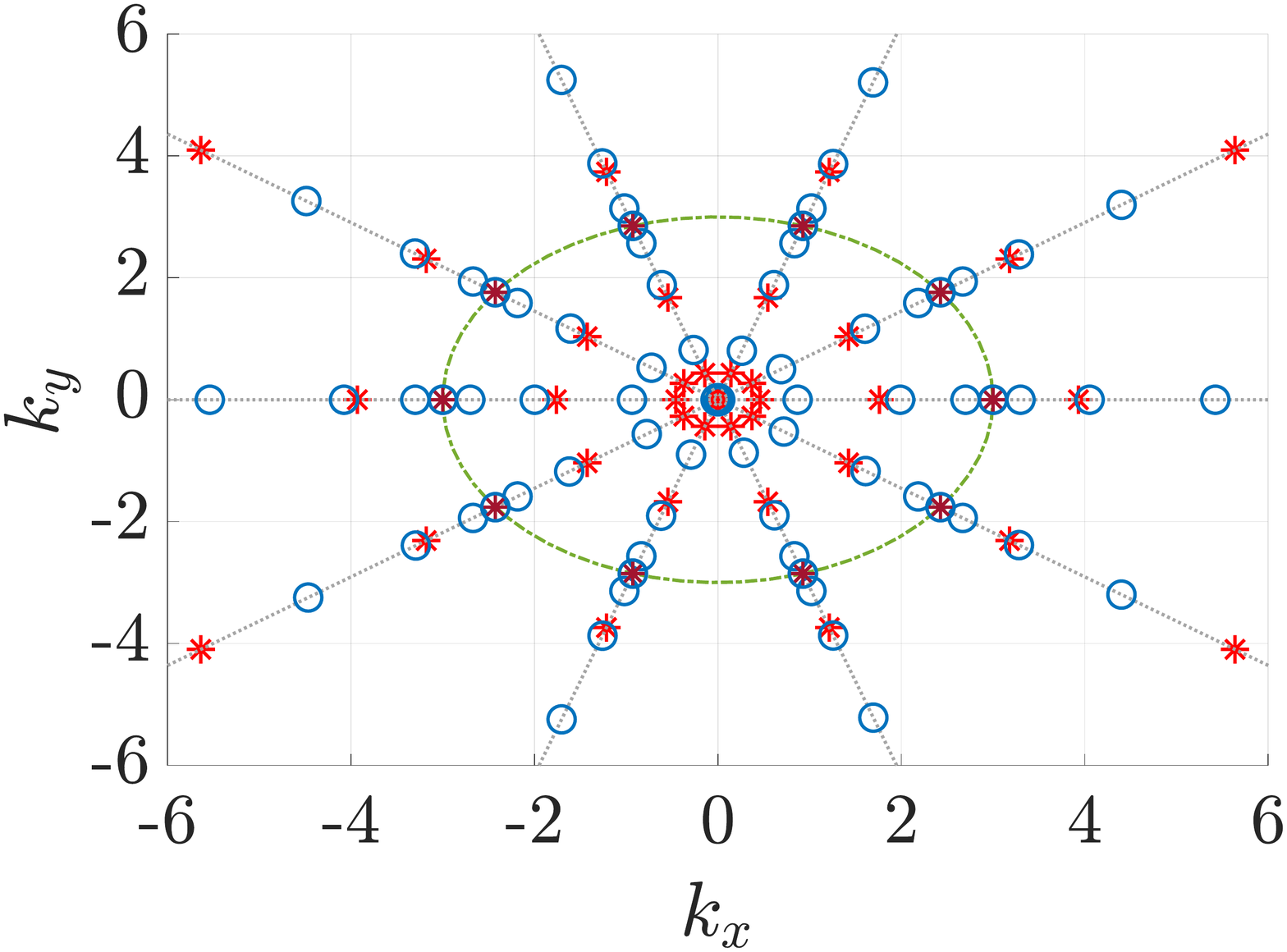}} \label{fig_interp_precomp_a} }%
    \hfill
    \subfloat[3d plot as panel (c).]{{\includegraphics[trim=1.5cm 0.0cm 3.0cm 1.5cm,width=0.45\textwidth]{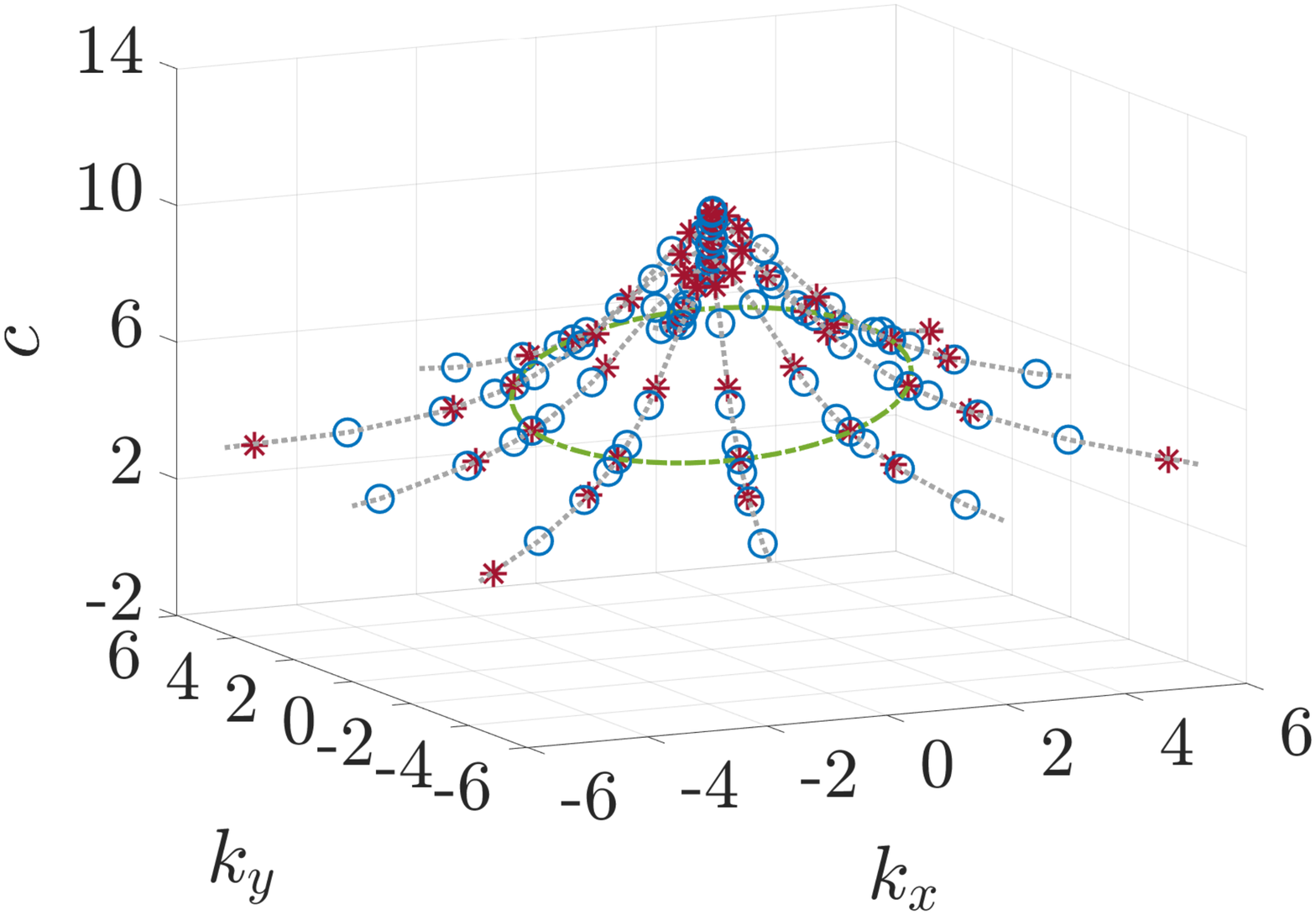} } \label{fig_interp_precomp_b}}\\%
    \subfloat[Planar plot showing interpolation of query point. Magenta crosses indicate the $k_q$ radius on each radial slice at angles $\sidx{\theta}{l}$, $\sidx{\theta}{l+1}$. Red astrisk indicates the query point $( k_q, \theta_q )$.]{{\includegraphics[trim=1.5cm 0.0cm 3.0cm 1.5cm,width=0.45\textwidth]{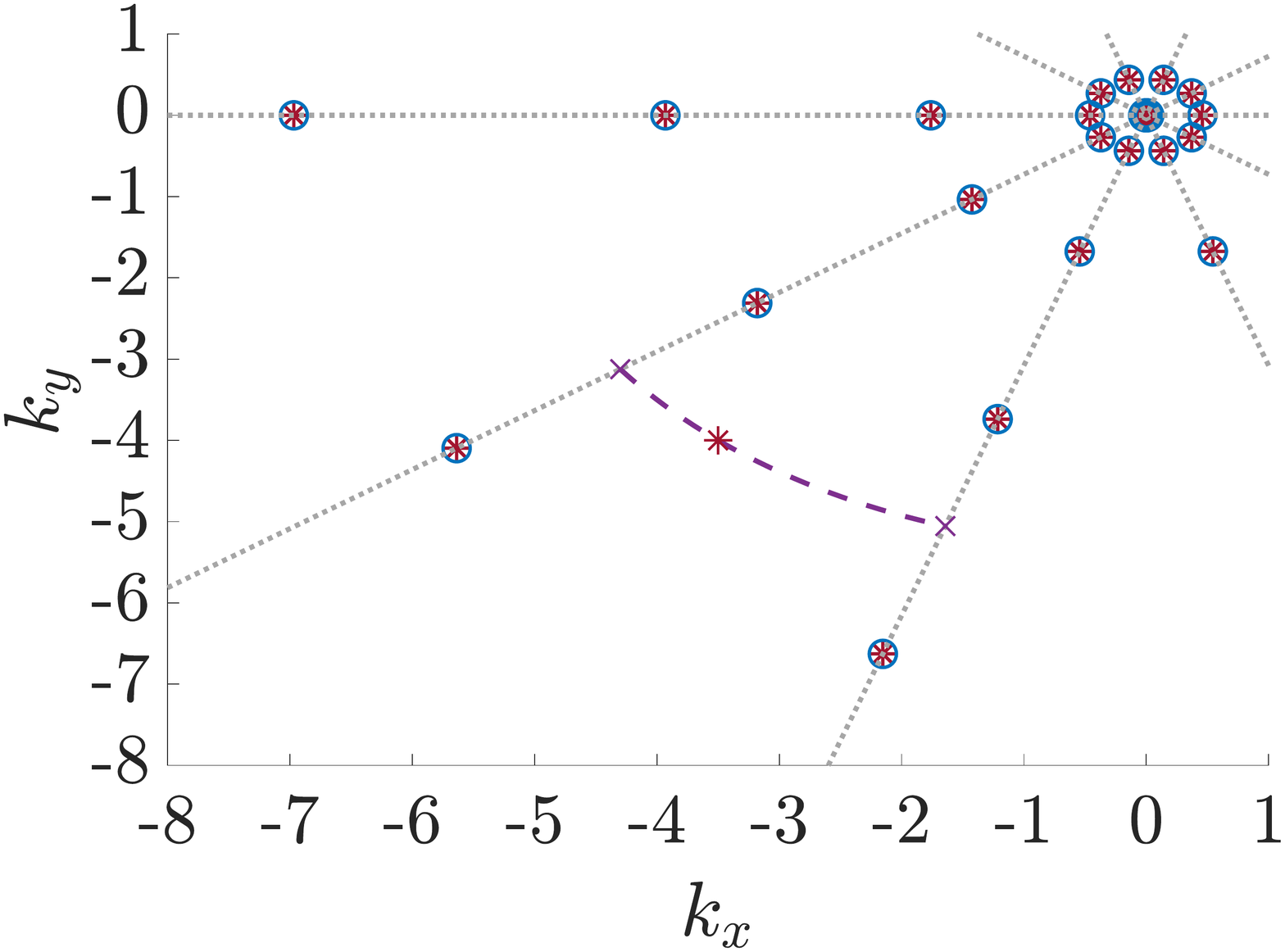}} \label{fig_interp_one_pt_a} }%
    \hfill
    \subfloat[3d plot as panel (e).]{{\includegraphics[trim=1.5cm 0.0cm 3.0cm 1.5cm,width=0.45\textwidth]{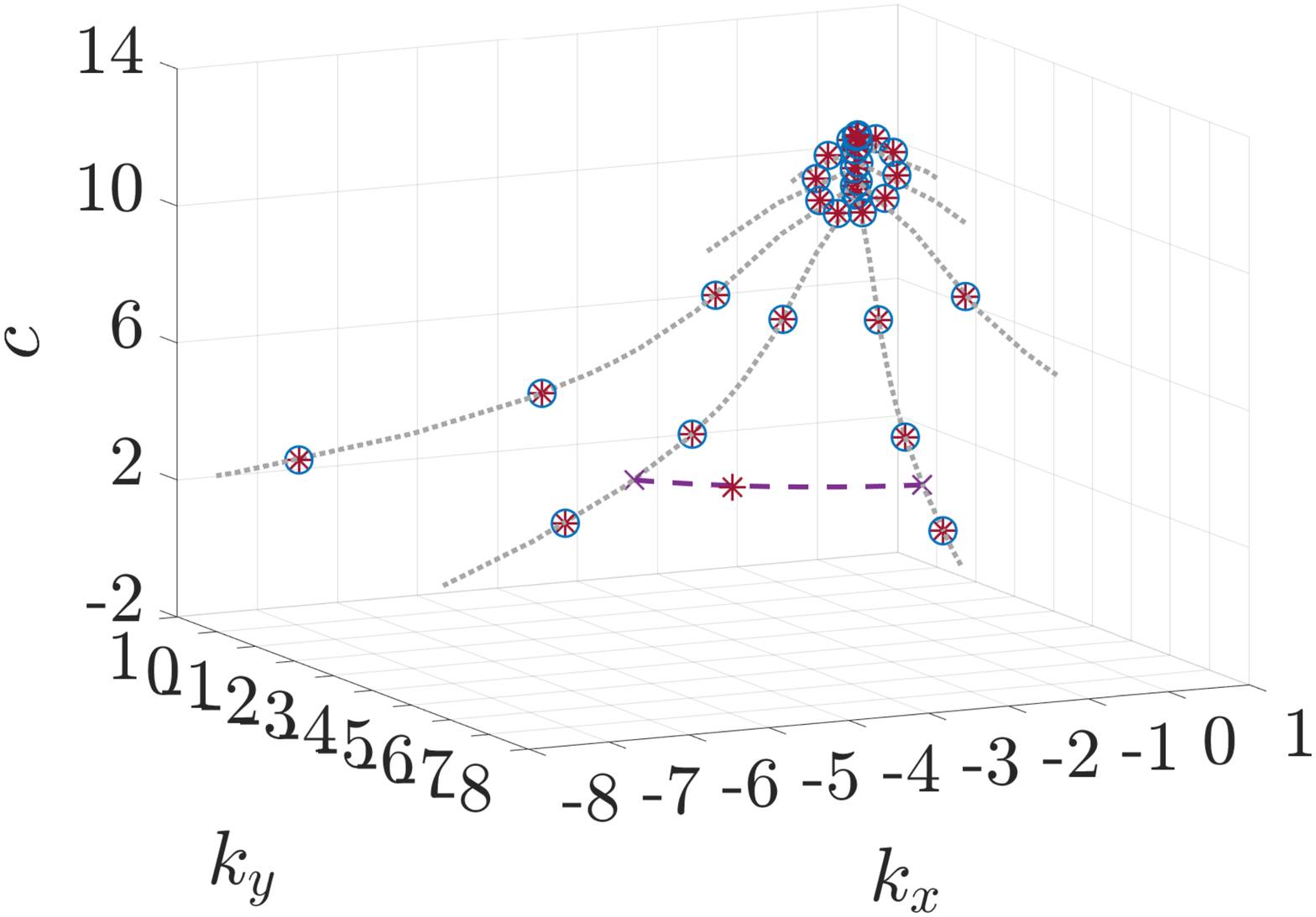} } \label{fig_interp_one_pt_b}}%
    \caption{Steps of PF-G.}%
    \label{fig_interp_scatter}%
\end{figure}%

\begin{pfgsteplist}
\item First, PF-R-a is executed at some nominal $k=\elminitval{k}$ and interpolation points at angles $\{ \sidx{\theta}{j} \}_{j=1}^{J}$ are calculated. See \cref{fig_interp_precomp_ang_a,fig_interp_precomp_ang_b}.\label{enum_step_i}
\item The results from \cref{enum_step_i} are used as the initial $\lavcinitval{v}$ values for PF-R-r calculating radially along each $\sidx{\theta}{j}$. The curve on each radial slice is then interpolated at predefined $\{ \sidx{k}{i} \}_{i=1}^{I}$ points. The control points for each radial slice are then replaced with the control points at these fixed $\sidx{k}{i}$ (we do not calculate new midpoint values). So there is now a 2-dimesional polar grid at angles $\sidx{\theta}{j}$ and radii $\sidx{k}{i}$. See arrangement in \cref{fig_interp_precomp_a,fig_interp_precomp_b}.\label{enum_step_ii}
\item For an arbitrary query point $( k_q, \theta_q )$ the nearest angles $\sidx{\theta}{l}$, $\sidx{\theta}{l+1}$ and radii $\sidx{k}{m}$, $\sidx{k}{m+1}$ are identified. The interpolant on radial slices at angles $\sidx{\theta}{l}$, $\sidx{\theta}{l+1}$ are calculated at radius $k_q$. \cref{eqn_PF_angular_c_system} is then used to calculate the angular derivatives. Finally, (cubic) interpolation is performed in an angular direction for angle $\theta_q$ to obtain the solution. See \cref{fig_interp_one_pt_a,fig_interp_one_pt_b}.\label{enum_step_iii}
\end{pfgsteplist}

Note that after \cref{enum_step_i,enum_step_ii} are calculated once, only \cref{enum_step_iii} need be performed for further query points, in a similar manner to PF-R-\{r,a\}.

There is a loss in accuracy because of the required use of cubic interpolation --due to not having the midpoint-- rather than the 4th order interpolation used in PF-R-\{r,a\}. However, for these purposes, it is not particularly significant. For clarity, we omit further analysis of PF-G: it is broadly similar to PF-R-r and does not add anything to the discussion.

\section{Convergence and error estimates}
\label{sec_accuracy}
It is well known that for sufficiently smooth solutions, spectral methods converge exponentially fast or with `spectral accuracy' \cite[ch.\ 1,2]{Boyd2001}. However, roundoff error poses a significant challenge for collocation methods due to the interaction of ill-conditioned matrices with commonly used \texttt{double} precision calculations \cite{Boyd2001}\cite{BaltenspergerTrummer2003}. We adopt a heuristic strategy to estimate the accuracy of each algorithm.

\subsection{Dependence of eigenvalue accuracy on order $N_z$}
\label{subsec_accuracy_eigenvalues_nz}
To determine accuracy depending on $N_z$, we first calculate a reference dispersion relation $R_\textrm{ref} = \{ ( \sidxlbl{k}{i}{ref}, \sidxlbl{c}{i}{ref} ) \}_{i=1}^{I}$ for the $\sidxlbl{k}{i}{ref}$ values distributed along the test interval $I_k$. This is done in high precision arithmetic, using the Advanpix library \cite{Advanpix}, for $N_\textrm{ref}=384$; this size of matrix exceeds what would be used in practice.

We calculate the relative normwise error in a candidate dispersion relation $R_\textrm{cand}$ (with $\sidxlbl{k}{i}{cand} = \sidxlbl{k}{i}{ref}$) as,
\begin{equation}
\epsilon := \frac{\norm{ [ \sbns \sidxlbl{c}{1}{cand} - \sidxlbl{c}{1}{ref}, \ldots, \sbns \sidxlbl{c}{J}{cand} - \sidxlbl{c}{J}{ref} \sbns ] }_{\infty}}{\norm{ [ \sbns \sidxlbl{c}{1}{ref}, \ldots, \sbns \sidxlbl{c}{J}{ref} \sbns ]   }_\infty}.\label{eqn_err_epsilon}
\end{equation}

\begin{figure}%
    \centering%
    \subfloat[Log-log plot of normwise relative error in candidate algorithms depending on $N_z$. The collocation method and path-following algorithm reduce error with spectral accuracy until roundoff error starts to dominate around $N_z \approx 65$. Path-following+MP which uses a high-precision initial result maintains the broadly the same low error despite the actual path-following calculations being preformed in \texttt{double} precision. DIM reduces error as predicted, as $\mathcal{O}(N_z^{-2})$.]{{\includegraphics[trim=1.5cm 0.0cm 3.0cm 1.5cm,width=0.45\textwidth]{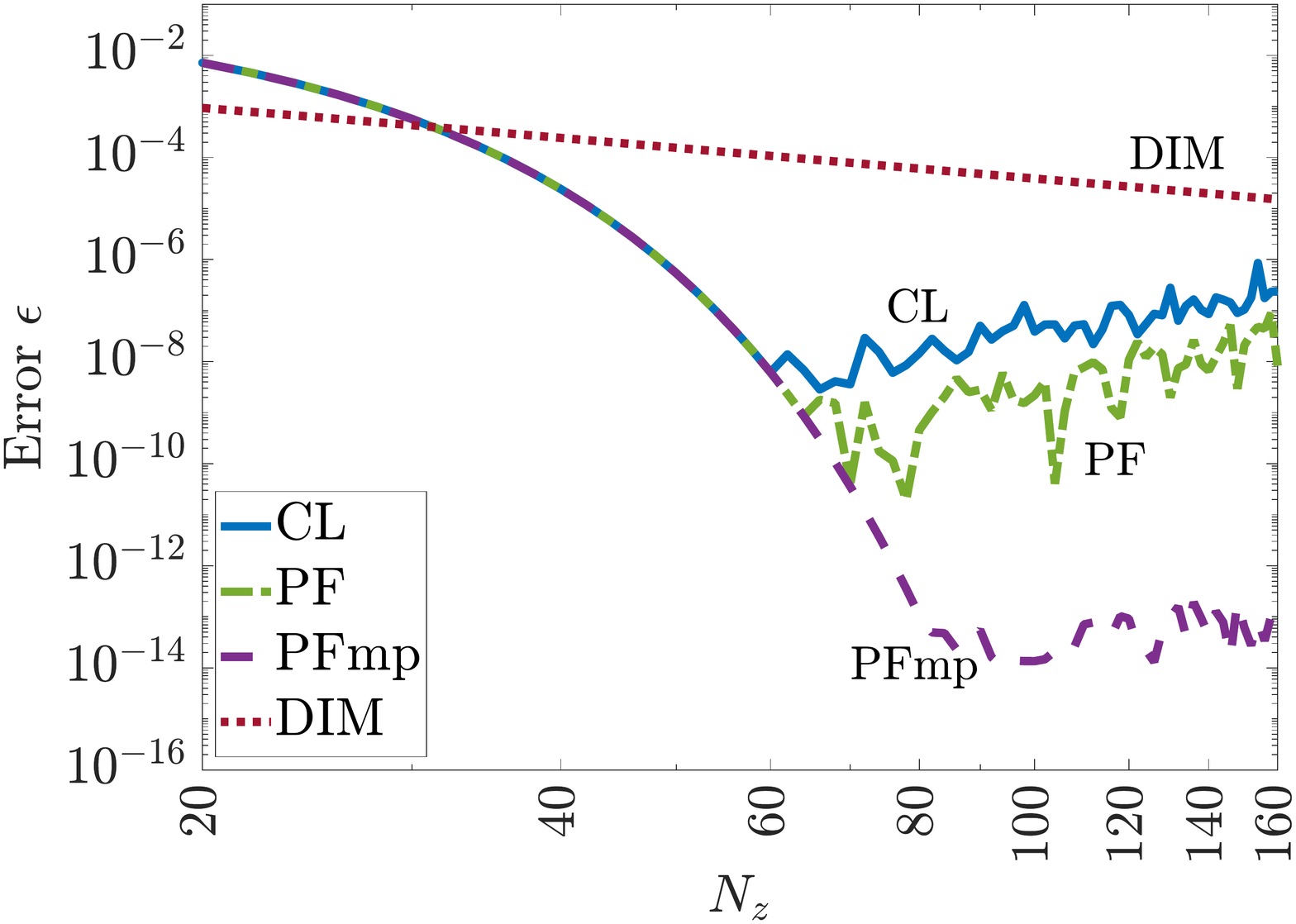}} \label{fig_eig_accuracy} }%
    \hfill
    \subfloat[Backwards error, condition number, and error estimate (backwards error $\times$ condition number). The time series for the quadratic eigenproblem solves for CL shown in blue, the linear solves for PF in magenta. From the condition numbers, shown with dotted lines, it can be seen that CL solves are slightly better conditioned but PF solves are of the same order of magnitude. The backwards error, shown in dashed lines, shows the linear solves in PF are appreciably more backwards stable. The error estimate clearly favours PF.]{{\includegraphics[trim=1.5cm 0.0cm 3.0cm 1.5cm,width=0.45\textwidth]{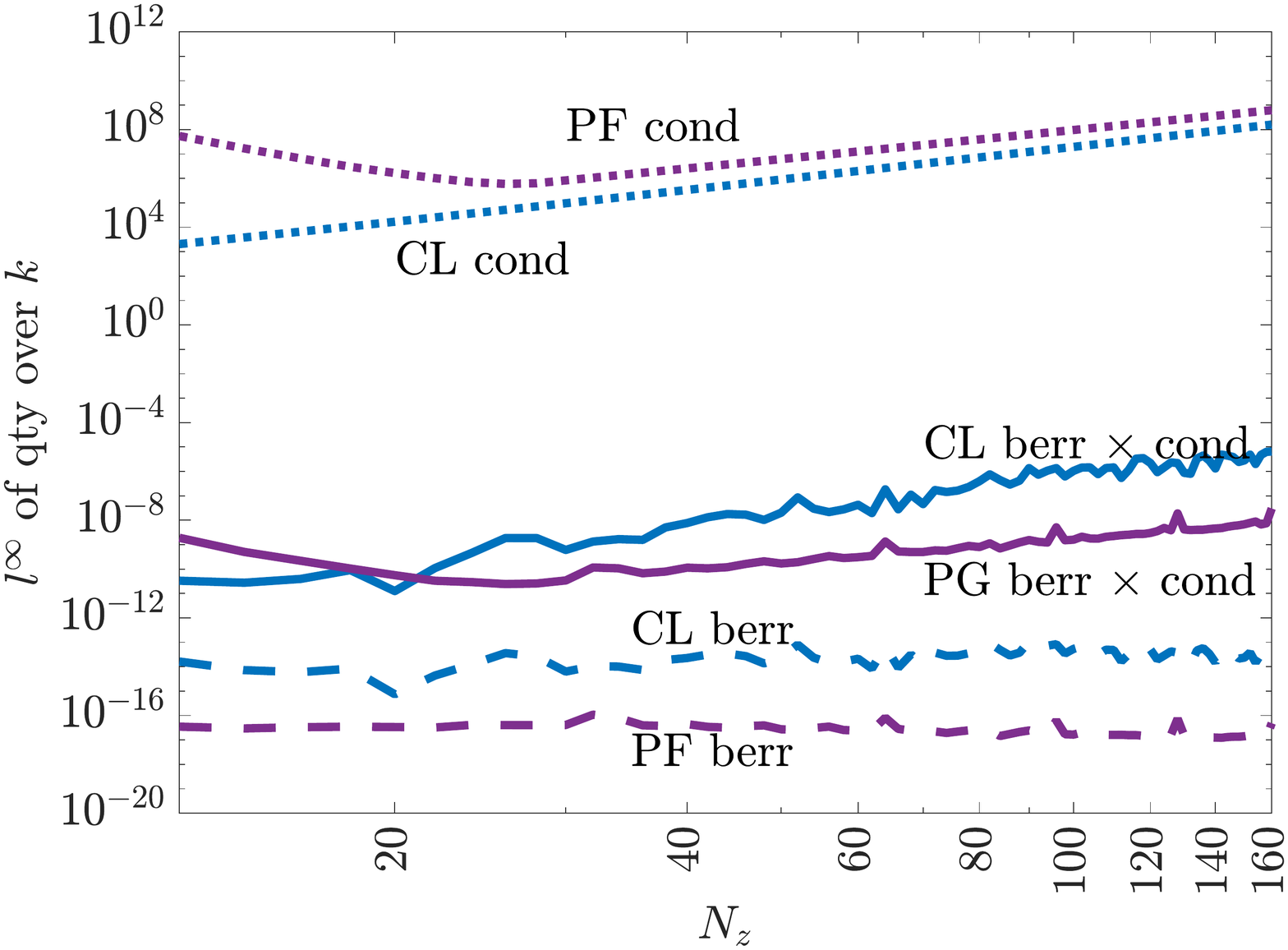} } \label{fig_berr_cond}}%
    \caption{Plots of error in eigenvalue computations and backwards error $+$ condition number estimates.}%
    \label{fig_acc_berr_cond}%
\end{figure}

This is done for for the $U_\textrm{T}$ shear profile using the CL-c, PF-R-r-c, PFmp-R-r-c, and DIM algorithms. The CL and PF methods reduce error with spectral accuracy until roundoff error starts to dominate. The high-precision initial calculation for the PFmp algorithm avoids this roundoff error and it can be seen that the path-following method itself retains this improved accuracy even in \texttt{double} precision.  DIM is included for indicative purposes.  See \cref{fig_eig_accuracy}.

A possible explanation for this can be found in comparison of the backwards error and conditioning of the quadratic eigenvalue solve used for the CL- schemes compared to the linear solves predominantly used in PF-. Although it is not a direct comparison ---the linear solves are used to calculate a derivative, not the value itself--- it may lend some insight. The backwards error for the linear solve can be calculated with \cite[eqn. 1.2]{HighamHigham1992} and the condition number in the usual manner:
\begin{equation}
\eta_{\textrm{L}} = \frac{\norm{ \lavc{b}-\lamx{M}\lamx{x}}_2 }{ \norm{\lamx{M}}_2 \norm{\lavc{v}}_2 + \norm{\lavc{b}}_2 }, \qquad \kappa_{\textrm{L}} = \norm{M^{-1}}_2 \norm{M}_2.
\end{equation}
The backwards error for the quadratic eigenproblem solve can be calculated using \cite[thm. 1, eqn. 2.3]{Tisseur2000} and the condition number using \cite[thm. 5, eqn. 2.15]{Tisseur2000}:
\begin{subequations}
\begin{gather}
\eta_{\textrm{Q}} = \frac{ \norm{ ( \lamxsb{A}{2} c^2 + \lamxsb{A}{1} c + \lamxsb{A}{0} ) \lavc{w} }_2 }{ ( \norm{\lamxsb{A}{2}}_2 |c|^2 + \norm{\lamxsb{A}{1}}_2 |c| + \norm{\lamxsb{A}{0}}_2 ) \norm{\lavc{w}}_2 },\\[0.5em]
\kappa_{\textrm{Q}} = \frac{ ( \norm{\lamxsb{A}{2}}_2 |c|^2 + \norm{\lamxsb{A}{1}}_2 |c| + \norm{\lamxsb{A}{0}}_2 ) \norm{\lavc{w}_l}_2 \norm{\lavc{w}}_2 }{|c||\lavc{w}_l^{*} ( 2\lamxsb{A}{2}c + \lamxsb{A}{1} ) \lavc{w} | }
\end{gather}
\end{subequations}
where $\lavc{w}_l$ is a corresponding left-eigenvector. Noting the usual inequality \cite[eqn. 1.3]{Tisseur2000},
\begin{equation}
\textrm{forward error } \le  \textrm{ condition number } \times \textrm{ backward error.}\label{eqn_berr_cond_ineq}
\end{equation}

For a range of $N_z$, we calculate the $\norm{\cdot}_{\infty}$ of the backwards errors $\eta_{\textrm{L}}$, $\eta_{\textrm{Q}}$  and condition numbers $\kappa_{\textrm{L}}$, $\kappa_{\textrm{Q}}$ over the $k$ vector. This then permits calculating the product from \cref{eqn_berr_cond_ineq}. This is shown in \cref{fig_berr_cond}. The condition numbers are of the same magnitude, $\kappa_{\textrm{L}} \approx \kappa_{\textrm{Q}}$, but the backwards error for the linear solves in PF-r is clearly smaller, $\eta_{\textrm{L}} \ll \eta_{\textrm{Q}}$. Although a direct comparison cannot be made, this suggests the path-following method has favourable numerical properties.

\subsection{$k$-dependent convergence}
\label{subsec_k_dep_cnvg}
\begin{figure}%
    \centering%
    \subfloat[Eigenvector plot for $k = \{ 0.025, \, 0.25, \, 2.5, \, 25, \, 250 \}$. As $k$ increases, the solution becomes numerically singular near $z=0$.]{{\includegraphics[trim=1.5cm 0.0cm 3.0cm 1.5cm,width=0.45\textwidth]{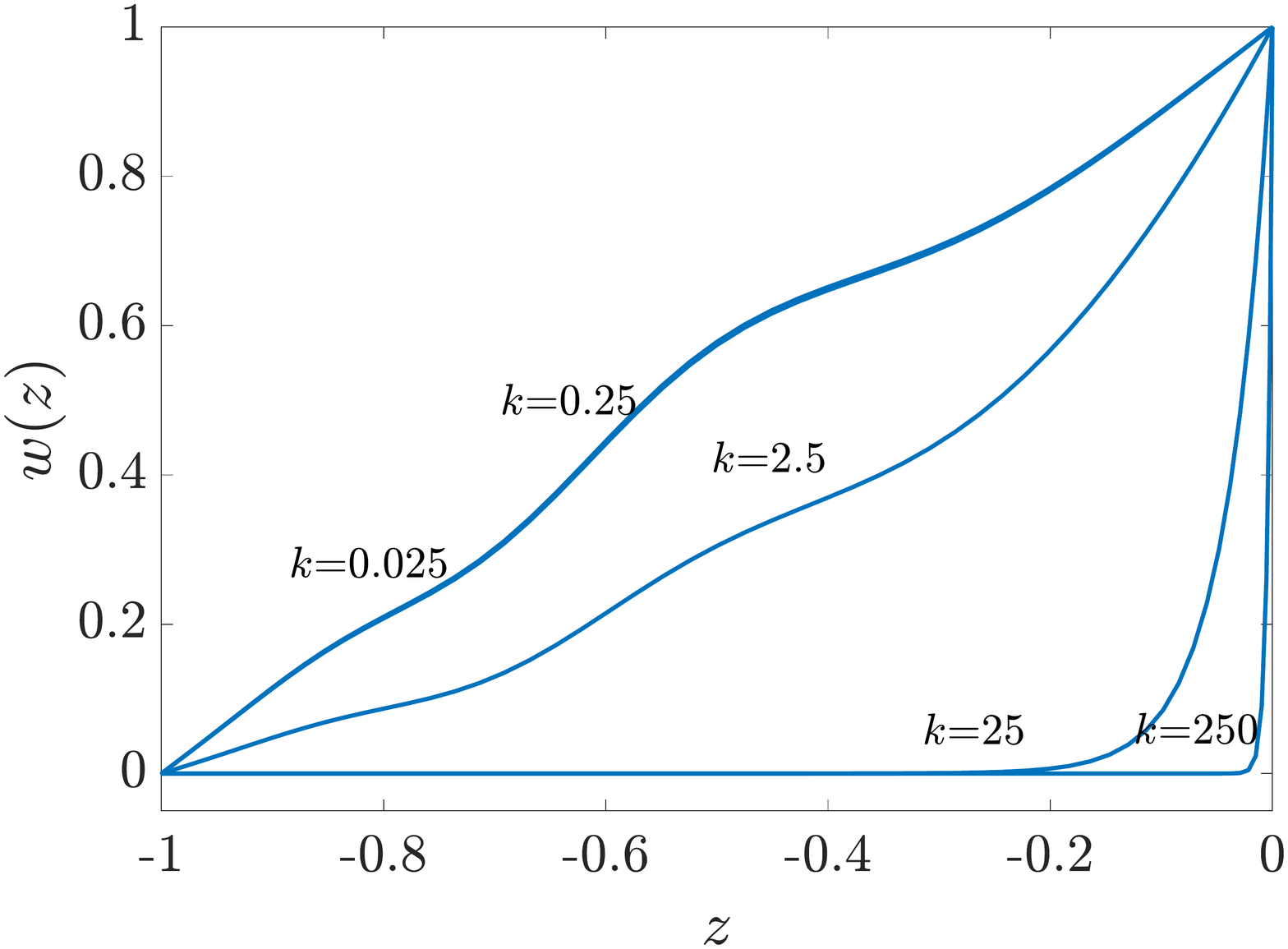}} \label{fig_evec_by_k} }%
    \hfill
    \subfloat[Lin-log plot showing required $N_z$ for Chebyshev series to converge depending on $k$. The increasingly numerically singular behaviour of the eigenvector requires much larger $N_z$ to reach convergence.]{{\includegraphics[trim=1.5cm 0.0cm 3.0cm 1.5cm,width=0.45\textwidth]{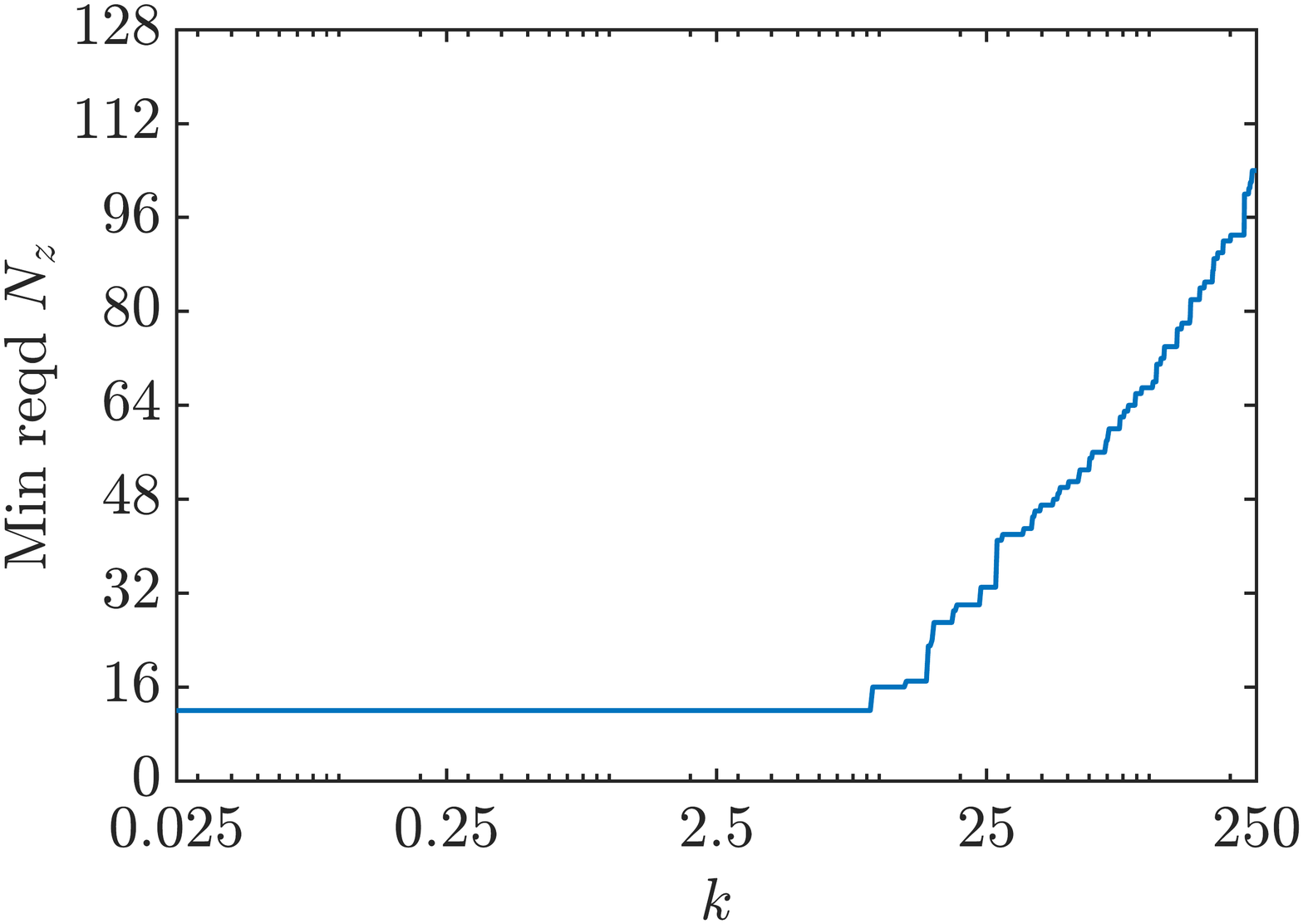} } \label{fig_k_convergence}}%
    \caption{Plot of eigenvectors for various $k$ and convergence properties.}%
    \label{fig_evec_cnvg}%
\end{figure}

As can be seen from \cref{fig_evec_by_k}, the eigenvectors become numerically singular at the surface as $k$ increases, implying that increasingly many basis polynomials are required to approximate the solution. This can be tested by using a similar algorithm as in \cite{AurentzTrefethen2017CCS} to determine when the Chebyshev series has converged. Specifically, we calculate an envelope then use a histogram to locate the plateau convergence region. We then determine the required $N_z$ to reach convergence for a range of $k$ values, as shown in \cref{fig_k_convergence}. For shorter wavelengths, much higher $N_z$ is required to reach convergence and so requiring more computational resources. This problem can be entirely ameliorated, as described in \cref{sec_pum}.

\section{Adaptive depth and partition of unity}
\label{sec_pum}
It is clear from the results in \cref{subsec_k_dep_cnvg} that as $k$ increases, the required $N_z$ becomes infeasibly large due to the singular behaviour of the eigenfunction. This can be avoided by using a smaller $h$ so that $h \ll 1$ for higher $k$, on the following rationale. We expect that the eigenfunction decays roughly as $e^{kz}$.  Therefore, we can estimate the depth below which the eigenfunction is effectively zero, for numerical purposes.  Let $\delta$ be the tolerance below which numerical values are considered zero, e.g.\ the ``machine epsilon''.  Let $h_{\delta}(k):= \min\{ 1, -\log(\delta) / k \}$ be an estimate of the depth, taking into account the finite depth, at which the eigenfunction decays below tolerance $\delta$ for a given $k$.

The CL-r scheme can be adapted for each calculate $k$. For a given $k$, we can set $h = h_{\delta}(k)$. The calculated eigenvalue for the phase velocity will be correct automatically. The eigenvector may be remapped back onto the original interval on any suitably large set of $z$ points chosen on the $[-1,0]$ interval using barycentric interpolation \cite{BerrutTrefethen2004}; the eigenvector will be zero for $z < -h_{\delta}(k)$.

This procedure becomes less obvious when considering the PF-r scheme because it would require remapping the entire system at each Runge--Kutta step.  To avoid this, we instead split the $k$ domain into several, partially overlapping, subintervals for which the chosen depth is suitable for all $k$ in that subinterval.  The path-following algorithm is then used on each subinterval independently with the appropriate depth.  To combine the subintervals and avoid loss of smoothness in the computed dispersion relation, a partition of unity method is used on the overlaps.

\begin{figure}%
    \centering%
    \includegraphics[trim=1cm 0cm 1.5cm 1cm,width=0.85\textwidth]{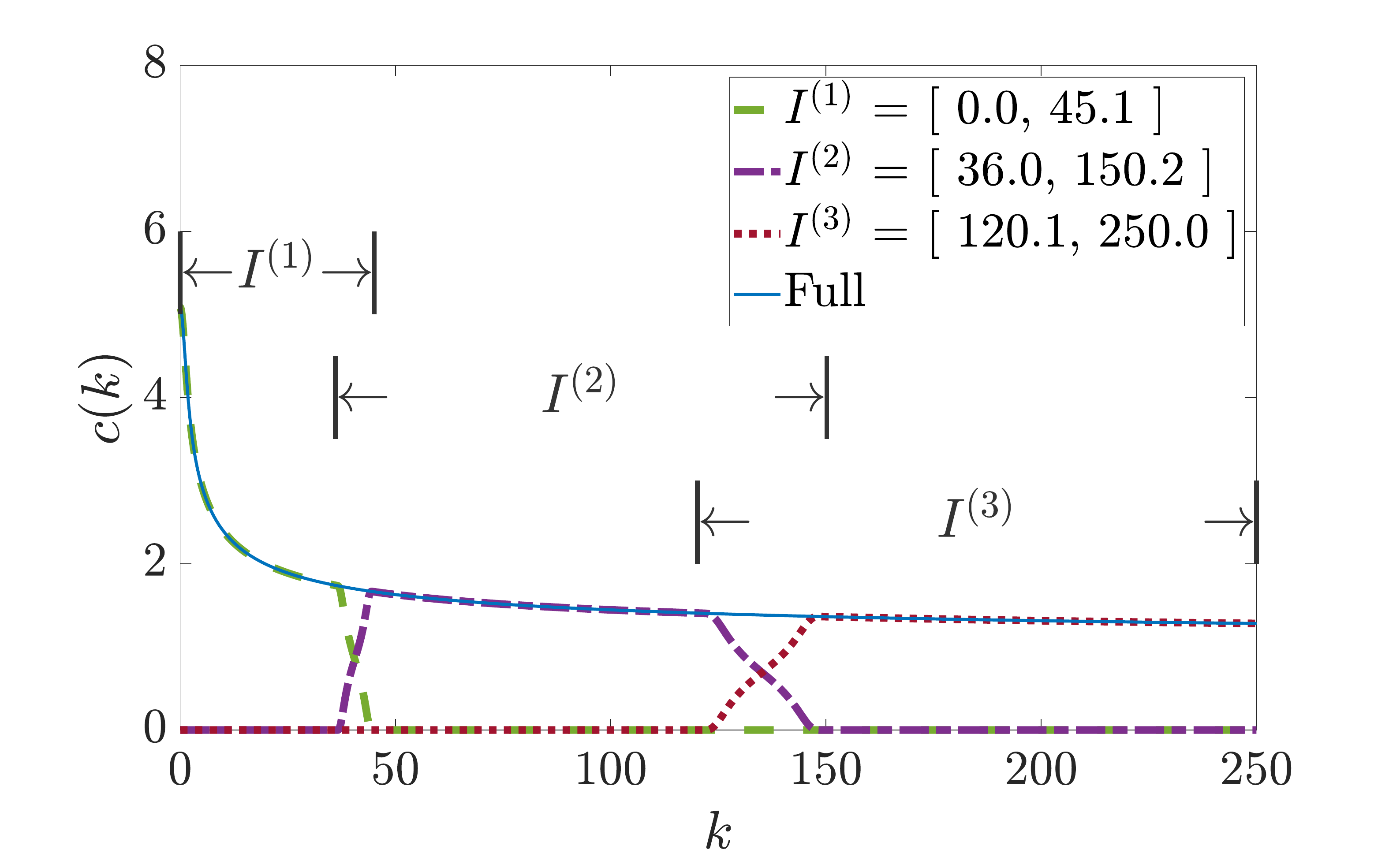}%
    \caption{Partition of unity method for dispersion relation.}
    \label{fig_pum_example}%
\end{figure}

We seek a scheme to choose subintervals $\sidx{I_k}{j} = [ \sidx{k_a}{j}, \sidx{k_b}{j} ]$ and corresponding depths $\sidx{h}{j}$ that is both simple and easy to implement. For some $\sidx{k}{j} \in \sidx{I_k}{j}$, we seek that $\sidx{h_{\min}}{j} \le  h_{\delta}(\sidx{k}{j}) \le \sidx{h_{\max}}{j}$ for $\sidx{h_{\min}}{j} = C_{\min} \sidx{h}{j}$ and $\sidx{h_{\max}}{j} = C_{\max} \sidx{h}{j}$ where $C_{\min}, C_{\max}$ are constants controlling the proportion of the $[\sidx{h}{j},0]$ interval that $h_{\delta}(\sidx{k}{j})$ should be within. In our computations, we found that $C_{\min} = 0.3$ and $C_{\max} = 0.8$ worked well. The subintervals and associated depths are then calculated as,
\begin{align*}
\sidx{I_k}{0} &= \left[ 0, \frac{\log( \delta )}{C_{\max}} \right], &&\sidx{h}{0} = 1,\\
\sidx{I_k}{j} &= \left[ \frac{\log( \delta )}{C_{\min}^{j-1}}, \frac{\log( \delta )}{\left( C_{\min}^{j} C_{\max} \right)} \right], &&\sidx{h}{j} = \min\left\{ 1, \frac{1}{2} \left( C_{\min}^{j-1} + C_{\min}^{j} C_{\max} \right) \right\}.
\end{align*}
This generates intervals $I_k^{(j)}$ in such a manner that $\sidx{k_b}{j} > \sidx{k_a}{j+1}$, i.e.\ there is some overlap in the intervals. We use the partition of unity method described in \cite{AitonDriscoll2018} to join the subintervals. This is demonstrated in \cref{fig_pum_example}.


\section{Performance analysis}
\label{sec_perf}
There are two variables which control the expected computation time for the candidate algorithms: the number of $z$ evaluation points, $N_z$, and the number of query points $N_q$. Since $N_z$ determines accuracy and is dependent on algorithm choice, we assume $N_z$ is set appropriately for each algorithm to achieve similar accuracy. Therefore, our primary concern shall be how the algorithms scale with $N_q$.

DIM will incur a fixed per-point computational cost that depends on the number of $z$ points, which we denote $\sigma_\textrm{DIM}(N_z)$. So, the expected cost is $\mathcal{O}( \sigma_\textrm{DIM}(N_z) N_q )$. Similarly, the collocation algorithm incurs a fixed per-point computational cost which also depends on the number of $z$ points used, $\sigma_\textrm{CL}(N_z)$. So, the expected cost is $\mathcal{O}( \sigma_\textrm{CL}(N_z) N_q )$. These estimates are valid for both the reduced and general problems.

In contrast, the reduced path-following algorithm incurs an initial computational cost dependent on the number of $z$ points, $\sigma_\textrm{PF-NI}(N_z)$, whereafter there is a very light-weight per-point cost, $\sigma_\textrm{PF-Q} \ll \sigma_\textrm{PF-NI}(N_z)$. Therefore, the expected computational costs is $\mathcal{O}( \sigma_\textrm{PF-NI}(N_z) + \sigma_\textrm{PF-Q} N_q )$. The general path-following algorithm is similar with only the coefficients changed. This is summarised in the following table, assuming the eigenvector output is not required:
\begin{center}
\begin{tabular}{|c|c|}
\hline
Algorithm & Computational Cost\\
\hline
DIM & $\mathcal{O}( \sigma_\textrm{DIM}(N_z) N_q )$ \\
CL-c \& CL-G-c & $\mathcal{O}( \sigma_\textrm{CL}(N_z) N_q )$ \\
PF-R-r-c &  $\mathcal{O}( \sigma_\textrm{PF-NI}(N_z) + \sigma_\textrm{PF-Q} N_q )$ \\
\hline
\end{tabular}
\end{center}
It immediately becomes clear that if $\sigma_\textrm{PF-NI}(N_z)$ is not too large and $\sigma_\textrm{PF-Q}$ is sufficiently small then as $N_q$ increases, the path-following algorithms are much more efficient.

\begin{figure}%
    \centering%
	\includegraphics[trim=0cm 0cm 2.0cm 0cm,width=0.98\textwidth]{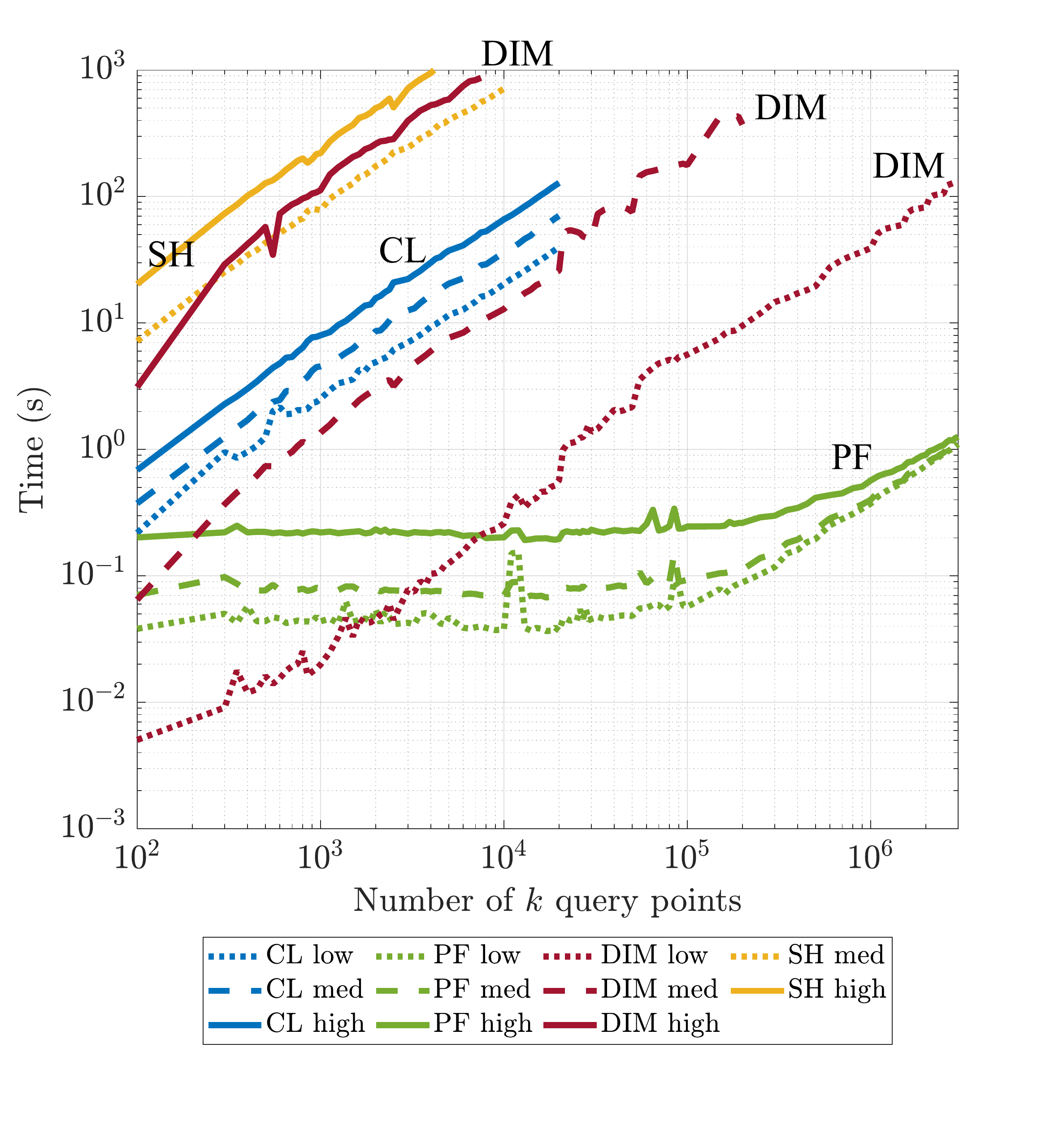}%
    \caption{Performance plot for reduced problem. Note that the collocation and DIM algorithms are clearly linear in complexity with respect to $N_q$. The path-following algorithm is also linear but this is only visible after around $10^{5}$ points because the interpolation cost is so small. For the $\approx 10^{-4}$ results, the path-following algorithm breaks-even at around 1200 query points; at $\approx 10^{-7}$ accuracy, it breaks-even at around 100 points; and, at $\approx 10^{-10}$ accuracy it is always faster. Asymptotically, the path-following algorithm is at least two orders of magnitude faster than both DIM and the collocation scheme.}%
    \label{fig_perf_std}%
\end{figure}

The asymptotic complexity claims are confirmed by practical testing. For clarity, we only test with the reduced problem. By measuring the time taken for each algorithm to compute the dispersion relation for differing $N_q$, we can determine the computational complexity in relation to $N_q$ as shown in in the log-log plot, \cref{fig_perf_std}. Each algorithm can be executed with different parameter choices that influence accuracy. As such, we calibrated each algorithm to produce output at three different accuracies (measured as relative normwise error using \cref{eqn_err_epsilon}): $\epsilon \approx 10^{-4}$, $\epsilon \approx 10^{-7}$, and $\epsilon \approx 10^{-10}$. As seen in the results, the path-following algorithm is asymptotically at least two orders of magnitude faster than both DIM and the collocation scheme.  The break-even point in $N_q$ at which the path-following scheme becomes faster than DIM is $N_q \approx 1200$ for $\epsilon \approx 10^{-4}$, $N_q \approx 100$ for $\epsilon \approx 10^{-7}$, and PF-R-r is always faster for $\epsilon \approx 10^{-10}$.

\section{Guidance on optimal parameter choices}
\label{sec_param}
Optimal parameter choices are predicated on two key properties: the required accuracy and the anticipated number of query points.

As can be observed from \cref{fig_perf_std}, the path-following algorithm is most effective when higher accuracy and at least a moderate number of query points are required. The nominal setup cost caused by the initial quadratic eigenproblem solve and numerical integration is dependent on the order of differentiation matrix used, $N_z$.  So this should be kept at the lowest value possible that maintains required accuracy. We found $N_z$ between 48 and 64 is optimal for the cases we tested. Furthermore, using $N_z$ too high risks roundoff error causing deleterious effects, cf.\ \cref{fig_eig_accuracy}.

The Dormand--Prince integrator requires a tolerance for the adaptive stepsize control.  We suggest that $10^{-11}$ is the smallest value to use when the initial eigenproblem solve is performed in \texttt{double} precision. If the initial eigenvalue solve can be performed more accurately, for example in high precision arithmetic, then the tolerance can be set around $10^{-15}$.  In any case, if using a smaller $N_z$ then the tolerance should be adjusted to match the accuracy from the collocation solution.

\section{Conclusions}
\label{sec_conclusions}
By considering the boundary value eigenproblem posed by the Rayleigh instability equation with linearised free-surface boundary condition as parameterised by wave number $k$, we can adapt the path-following scheme in \cite{LoiselMaxwell2018} to efficiently calculate the dispersion relation at high accuracy.  This efficiency is achieved by first exchanging many expensive QZ decompositions on a size $2N$ matrix for one QZ decomposition and some linear solves on a size $N$ matrix; secondly, we `front load' the computational cost into the numerical integration with light-weight Hermite interpolation being used to compute the sought solution points.

The accuracy tests in \cref{sec_accuracy} suggest that the path-following algorithm can maintain the accuracy of the initial eigenpair $\lavcinitval{v}$ and appears to be numerically more stable than the QZ decomposition used to obtain the initial eigenpair.

The algorithm is extended to permit calculation in the 2d $\vc{k}$-plane with scattered data and some preliminary discussion of critical layers is given.

In other tests, not included here, it is clear the same approach works well for other problems from physics and engineering, assuming the problem is parametrised by a real scalar.  Additional difficulties arise when there are exceptional points or bifurcations in the solution curve, or if the ODEs become stiff.  These challenges may form the basis of future work.

The MATLAB library used to perform the calculations used in this paper is maintained at \cite{NessieWaterWaves}.

\FloatBarrier
\bibliographystyle{siamplain}

\end{document}